\documentclass[3p,preprint]{elsarticle_preprint}

\usepackage{amsmath}
\usepackage{amsfonts}
\usepackage{amssymb,latexsym}
\usepackage{amsthm}
\usepackage{graphicx}

\newtheorem{definition}{Definition}[section]
\newtheorem{theorem}{Theorem}[section]

\usepackage{color}
\usepackage{newlfont}
\usepackage{subfigure}
\usepackage{verbatim}
\usepackage{rotating}
\usepackage{multirow}
\usepackage{units}
\usepackage{bm}
\usepackage[english]{babel}
\usepackage[utf8]{inputenc}
\usepackage{algorithm}
\usepackage[noend]{algpseudocode}
\usepackage{booktabs}
\usepackage{caption}

\journal{G. Garmanjani, R. Cavoretto, M. Esmaeilbeigi}

\begin{document}

\begin{frontmatter}

\title{A RBF partition of unity collocation method based on finite difference for initial-boundary value problems}


\author[address-MA]{G. Garmanjani}
\ead{gholamreza.garmanjani@stu.malayeru.ac.ir, reza.garmanjani@gmail.com}

\author[address-TO]{R. Cavoretto}
\ead{roberto.cavoretto@unito.it}

\author[address-MA]{M. Esmaeilbeigi\corref{corrauthor}}
\cortext[corrauthor]{Corresponding author}
\ead{m.esmaeilbeigi@malayeru.ac.ir, esmaeelbeigi@yahoo.com}

\address[address-MA]{Department of Mathematics, Faculty of Mathematics Science and Statistics, Malayer University, Malayer 65719-95863, Iran}
\address[address-TO]{Department of Mathematics \lq\lq Giuseppe Peano\rq\rq, University of Torino, via Carlo Alberto 10, 10123 Torino, Italy}

\begin{abstract}
Meshfree radial basis function (RBF) methods are popular tools used to numerically solve partial differential equations (PDEs). They take advantage of being flexible with respect to geometry, easy to implement in higher dimensions, and can also provide high order convergence. Since one of the main disadvantages of global RBF-based methods is generally the computational cost associated with the solution of large linear systems, in this paper we focus on a localizing RBF partition of unity method (RBF-PUM) based on a finite difference (FD) scheme. Specifically, we propose a new RBF-PUM-FD collocation method, which can successfully be applied to solve time-dependent PDEs. This approach allows to significantly decrease ill-conditioning of traditional RBF-based methods. Moreover, the RBF-PUM-FD scheme results in a sparse matrix system, reducing the computational effort but maintaining at the same time a high level of accuracy. Numerical experiments show performances of our collocation scheme on two benchmark problems, involving unsteady convection-diffusion and pseudo-parabolic equations.
\end{abstract}


\end{frontmatter}



\section{Introduction}\label{sec1}

A number of partial differential equations (PDEs) arises as initial-boundary value problems in various scientific and engineering applications. The latter include processes with slow evolution as heat conduction and diffusion, changes in time of processes in life and social sciences, atmospheric processes, ocean circulation, flow in wind tunnels, lee waves, eddies, etc. Analytic and numerical treatments of this class of problems involve several complications due to possible non-linearities of the studied processes, instabilities of numerical schemes causing blow-ups in simulations, and increases of the amount of data due to extra dimensionality. Hence, from theoretical and applied angles, remaining open questions have released active investigations of such models in mathematical and numerical modeling, as well as in large-scale scientific computing. Also, since most of the PDEs cannot be solved exactly, developing accurate and efficient approximation schemes for computing numerically the solution of differential equations has actively been required in several research areas.

In this study, we focus on two difficult benchmark problems, including an unsteady convection-diffusion equation and a pseudo-parabolic problem.

The convection-diffusion equation assumes a very important role in fluid flow and heat transfer problems. In fact, this equation is peculiar because it is a combination of two dissimilar phenomena as convection and diffusion, where many fluid variables such as mass, heat, energy and vorticity can be involved. Additionally, it can also be viewed as a simplified model problem to the Navier-Stokes equations, thus making the numerical prediction of the solution of the convection-diffusion equation very significant in computational fluid dynamics \cite{cha07}. The convection-diffusion equation has been then extensively used to model diverse processes in science and engineering, such as pollutant dispersal in a river estuary, atmospheric pollution simulation, ground water transportation, water transfer in soil, fluid motion, heat transfer, astrophysics, oceanography, meteorology, semiconductors, hydraulics, pollutant and sediment transport, chemical engineering, and the intrusion of salt water into fresh water aquifers, to name a few. For more details, the reader can refer to \cite{bea60,ise72,mor96,par80, roa76, pat80, bea72}. In such references, mathematical properties of convection-diffusion equation have been widely checked and many discussions have been done. In the literature, excellent summaries have been provided in \cite{mor96,baz08,el88,eno92,eva85,tem09}. Among the numerical methods used in simulation to solve convection-diffusion problems, we can mention different approximating schemes such as the finite difference (FD) method \cite{gol10,sun14}, the finite volume (FV) scheme \cite{cou10,cla13}, the finite element (FE) technique \cite{coc98,zhu11}, and the meshfree radial basis function (RBF) method \cite{ahm17,cha07,chi06}.

The linear pseudo-parabolic equation was first examined by S.L. Sobolev \cite{sob54} in 1954; for this reason, it is also named as a Sobolev type equation. The pseudo-parabolic equation arises in the flow of fluids through fissured rock \cite{bar60}, thermodynamics \cite{che68} shear in second order fluids \cite{tin63}, soil mechanics \cite{tay52}, and other applications. Such equations are important examples of PDEs, which contain a third order mixed derivative with respect to time and space. They are used to characterize wave motion, which are momentous not only in hydrodynamics but also in many other areas of engineering and science \cite{gao17}. From a theoretical viewpoint several results are known about existence, uniqueness, and properties of solutions. Brill \cite{bri77} and Showalter \cite{sho72} established the existence of solutions of semilinear evolution equations of Sobolev type in Banach space, while the global solvability and blow-up of equations of Sobolev type were considered in \cite{als11}. Numerically, a lot of simulation methods have been developed for pseudo-parabolic equations \cite{zha16,sol14}, but most of the existing works are based on the classical FE methods\cite{shi16,gao17}, FD schemes \cite{sun02}, or FV element methods \cite{yan08} as discretization tools.

As computer capabilities kept improving throughout the last three decades, it became possible to solve more and more complicated problems. The enhancement of possibilities to solve such complex physical end engineering problems has also been due to advances in numerical methods and improvements of algorithms in terms of efficiency. For instance, it became possible to simulate large-scale problems. Most of engineering and physical problems are solved by FD schemes, FE techniques, control volume (CV) or boundary element (BE) methods. All these numerical methods are based on a grid discretization that has to be generated in advance or dynamically improved as the solution progresses through adaptive meshing. Nevertheless, although it is widely known that mesh generation remains one of the biggest challenges in mesh-based methods \cite{fass07}, sometimes even the most computationally severe problems can be solved accurately provided that an acceptable mesh is found. So, even now, the production of appropriate mesh for such complex problems remains the main task and can take most of the computational time and effort needed to determine the approximated solution of the PDEs. Furthermore, not only mesh-based methods can be extremely complicated, but for many problems or dynamic processes such as those with large deformations they are known not to be cost-effective.

Over the recent years, significant effort has been dedicated to the study of so-called meshfree methods, which are also cited in the literature as meshless, element-free, gridless or cloud methods. The purpose of meshless methods is to eliminate at least the structure of the mesh and approximate the solution entirely using the nodes/points as a scattered or quasi-random set of points rather than nodes of an element/grid-based discretization \cite{zer00}. A class of meshless methods are based on RBF collocation methods which put radial functions as the basis functions for the collocation \cite{kan90a,kan90b}.

Recently, the application of meshless RBF-based methods has gained popularity in science community for a number of reasons. The most prominent characteristics of these methods are meshfree character, flexibility in dealing with complex geometries, spectral convergence property, and easy extension to multi-dimensional problems. In fact, a RBF only depends upon the distances between centers, so that RBF interpolants are typically used to approximate given functions. Additionally, the RBFs are also used as a kernel in support vector classification. However the use of globally supported basis functions leads to fully-populated collocation matrices, which become increasingly ill-conditioned and computationally expensive with increasing data set size. These constraints have thus motivated researchers and scientists to investigate different methods for restricting the influence region of the basis functions, thereby mitigating the computational cost and numerical conditioning issues, while maintaining the efficiency and flexibility of the previous formulation. Relating to this field several excellent books were written by different authors, who analyzed various properties and possible uses of RBFs \cite{buh03,fass07,fass15,wen05}.

The simplest way of making a RBF collocation method localizing is to restrict the support of the basis functions. In this case the standard globally supported basis functions are replaced by functions which are nonzero only within their support. In this way, the resulting global collocation matrices are sparsely populated, and may be solved far more efficiently than their fully-populated global counterparts. Wendland \cite{wen95} and Wu \cite{wu95} constructed compactly supported schemes, proposing a series of positive definite polynomial functions for various spatial dimensions and orders of continuity. Such methods are straightforward to implement, but they are not much flexible and it can be shown that large support sizes are often required in order to achieve reasonable accuracy and convergence rates \cite{flo96}. An alternative approach is to perform a direct domain decomposition of the RBF collocation, whereby the solution domain is divided into several interacting subdomains which may or may not overlap depending on the formulation. While domain decomposition is attractive in principle, the definition of suitable subdomains may be quite complex and highly problem-dependent \cite{ste13}.

Here, we now deal with the two main techniques currently used to reduce the computational cost of RBF-based methods through localization. The first approach, named RBF-FD method, can be seen as a generalization of FD schemes (but possibly with stencils supported on scattered point sets) and is given by a combination of RBFs and FD. To the best of our knowledge, Tolstykh \cite{tol00} was the first to publish the method, which then has been widely studied, see e.g. \cite{cha07,for15}. The second approach, which is the focus of this paper, is the RBF partition of unity method (RBF-PUM). It -- also known as RBF-PU method -- is performed by decomposing the problem domain in a suitable number of overlapping subdomains or patches forming a covering of the original domain. Then, after constructing a local RBF approximant on each subdomain, such approximations are weighted together by compactly supported PU weight functions to form the global fit. With this scheme a large problem is decomposed into many small problems and therefore we can work with a large number of points. Also, the convergence properties of the local approximations can be leveraged, while local couplings between approximations on different patches are enforced through the PU framework. PU methods for solving PDEs were originally proposed by Babu$\check{\text{s}}$ka and Melenk \cite{mel96,bab97}, suggesting the idea of combining RBF approximations with partition of unity. Later the PUM combined with compactly supported RBFs (CSRBFs) was introduced and analyzed by Wendland \cite{wen02} for interpolation purposes, while more recently Cavoretto, De Rossi et al. have proposed several efficient algorithms for RBF-PUM interpolation of scattered data \cite{cav15a,cav16,cav17a,cav17b,cav17c}. Finally, Larsson et al. have published some papers about RBF-PU collocation method for PDE problems \cite{saf15,shc16}, also providing a few theoretical results for the approximation errors of such an approach.

Generally, various FD schemes can be used for the time discretization in solving initial-boundary value problems. These schemes can be categorized into three classes: explicit, implicit, and semi-implicit schemes. The explicit scheme is efficient but has poor stability, making its applicability very difficult with a large time step; consequently, the overall computational cost is usually high. On the other hand, the implicit scheme is unconditionally stable and allows a larger time step, but requires to solve a linear algebraic system at any time step. The semi-implicit scheme is somewhere between the explicit and implicit schemes. Although the conclusion as to which type of scheme should be chosen is to consider case by case, usually the implicit scheme is preferred for stiff initial value problems \cite{lia12,esm16}.

In this paper, we are interested in solving time-dependent PDEs in a two-dimensional space. Specifically, we propose a new method, called RBF-PUM-FD, that is a local RBF-PU collocation method based on a FD scheme for solving convection-diffusion and pseudo-parabolic problems. The advantage of this approach is that -- due to use of a FD scheme -- the linear or nonlinear PDEs results in a system of linear equations, whereas -- thanks to application of PUM -- a large problem is decomposed into many small problems. Therefore, in the approximation procedure we can also deal with a large number of points. This method is proposed to cope with ill-conditioning, overcoming the traditional issues of global RBF-based methods. Moreover, using an implicit FD technique for the time discretization leads to an unconditionally stable method, which requires the solution of a sparse linear system. In this work, we thus have a two-fold scope: first, exhibiting efficiency of our RBF-PUM-FD method in terms of CPU times; then, analyzing its numerical accuracy and stability. Finally, in order to show goodness of our approach, we compare our collocation method with two existing computational techniques: the global RBF-FD method proposed in \cite{esm16} and the RBF-PUM studied in \cite{saf15}.

The layout of the article is as follows. In Section \ref{sec2}, we review the main theoretical results concerning RBF approximation. Then we briefly present the PUM collocation by RBFs in next section. In Section \ref{cd_pp_prbs}, we propose the RBF-PUM collocation method based on FD for time-dependent PDEs. In Section \ref{sec5}, we provide a stability analysis of our collocation method. The results of our extended numerical experiments are presented in Section \ref{sec6}. Finally, Section \ref{sec7} is dedicated to brief conclusions and future work.


\section{RBF approximation}\label{sec2}

Since the numerical solution of time-dependent PDEs by RBF methods is based on a scattered data interpolation problem, in this section we review the main theoretical aspects concerning interpolation (or collocation) via RBFs. The associated schemes are meshfree and effectively works also with scattered or irregularly distributed data points, therefore they turn out to be flexible in terms of the geometry of the domain \cite{buh03}.


\subsection{RBF collocation by conditionally positive definite functions} \label{sec21}

In scattered data interpolation using RBFs the approximation of a function $u(\boldsymbol{x})$ at distinct \textsl{data points} or \textsl{nodes} or \textsl{centers} $X=\lbrace \boldsymbol{x}_1,\ldots,\boldsymbol{x}_N \rbrace$, defined in a domain $\Omega \subseteq \mathbb{R}^d$, may be written as a linear combination of RBFs, also adding certain polynomials to its expansion. The corresponding RBF interpolation problem consists thus in finding an interpolant $s_{u,X}: \Omega\rightarrow \mathbb{R}$ of the form
\begin{equation} \label{eq:21} 
s_{u,X}(\boldsymbol{x})=\sum_{j=1}^N \alpha_j \phi_{\epsilon}(||\boldsymbol{x}-\boldsymbol{x}_j||)+\sum_{k=1}^Q\beta_{k} p_{k}(\boldsymbol{x}), \qquad \boldsymbol{x} \in \mathbb{R}^d,
\end{equation}
where $\alpha_j$ and $\beta_k$ are unknown real coefficients, $||\cdot||$ denotes the Euclidean distance, and $\phi_{\epsilon}: \mathbb{R}_{\geq 0} \rightarrow \mathbb{R}$ is a RBF depending on a \textsl{shape parameter} $\epsilon > 0$ such that $\phi_{\epsilon}(||\boldsymbol{x}-\boldsymbol{z}||)=\phi(\epsilon ||\boldsymbol{x}-\boldsymbol{z}||)$, for all $\boldsymbol{x},\boldsymbol{z} \in \Omega$. For simplicity, from now on we refer to $\phi_{\epsilon}$ as $\phi$. In Table~\ref{tab_rbf} we report a list of some well-known RBFs with their orders of smoothness \cite{fass07}. 

\begin{table}
\begin{center}
\begin{tabular}{ll}
\hline
\rule[0mm]{0mm}{3ex}
RBF  & $\phi_{\epsilon}(r)$ \\
\hline
\rule[0mm]{0mm}{3ex}
{Gaussian $C^{\infty}$} (GA) & ${e}^{-\epsilon^2 r^2}$   \\
\rule[0mm]{0mm}{3ex}
{MultiQuadric $C^{\infty}$} (MQ) & $(1+\epsilon^2r^2)^{1/2}$   \\
\rule[0mm]{0mm}{3ex}
{Inverse MultiQuadric $C^{\infty}$} (IMQ) & $(1+\epsilon^2r^2)^{-1/2}$   \\
\rule[0mm]{0mm}{3ex}
{Thin Plate Spline $C^{\nu + 1}$} (TPS) & $(-1)^{\nu+1}r^{2\nu}\log r$   \\
\rule[0mm]{0mm}{3ex}
{Mat$\acute{\text{e}}$rn $C^4$} (M4)  & ${e}^{-\epsilon r} (\epsilon^2r^2+3\epsilon r+3)$   \\
\rule[0mm]{0mm}{3ex}
{Mat$\acute{\text{e}}$rn $C^2$} (M2) & ${e}^{-\epsilon r} (\epsilon r+1)$   \\
\rule[0mm]{0mm}{3ex}
{Wendland $C^4$} (W4)  & $\left(1-\epsilon r\right)_+^6(35\epsilon^2 r^2+18\epsilon r+3)$  \\
\rule[0mm]{0mm}{3ex}
{Wendland $C^2$} (W2) & $\left(1-\epsilon r\right)_+^4\left(4\epsilon r+1\right)$  \\
\hline
\end{tabular}
\end{center}
\caption{Examples of some popular RBFs, where $r=||\cdot||$ is the Euclidean norm, $(\cdot)_+$ denotes the truncated power function, and $\nu \in \mathbb{N}$.}
\label{tab_rbf}
\end{table}

Moreover, in Eq. \eqref{eq:21} $p_1,\ldots,p_Q$ form a basis for the $Q$-dimensional space $\pi_{m-1}(\mathbb{R}^d)$ of polynomials of total degree $\leq m-1$ in $d$ variables. To cope with additional degrees of freedom, the interpolation conditions
\begin{equation*} 
s_{u,X}(\boldsymbol{x}_i)=u(\boldsymbol{x}_i), \qquad  i = 1,\ldots, N,
\end{equation*} 
are completed by the additional conditions
\begin{equation} \label{eq:23} 
\sum^N_{j=1}{{\alpha }_j}p_k(\boldsymbol{x}_j)=0, \qquad k = 1,\ldots, Q. 
\end{equation} 
Solving the interpolation problem based on the expansion \eqref{eq:21} amounts therefore to solving the system of linear equations
\begin{equation} \label{eq:24} 
\left( \begin{array}{cc}
A & P \\ 
P^T & O \end{array}
\right)\left( \begin{array}{c}
\boldsymbol{\alpha}  \\ 
\boldsymbol{\beta}  \end{array}
\right)=\left( \begin{array}{c}
\boldsymbol{u} \\ 
\boldsymbol{0} \end{array}
\right),                                                                                     
\end{equation} 
where the matrix $A\in {{\mathbb R}}^{N\times N}$ is given by $A_{ij}=\phi(||\boldsymbol{x}_i-\boldsymbol{x}_j||)$, $i,j=1,\ldots,N$, $P\in {{\mathbb R}}^{N\times Q}$ has entries $P_{jk}=(p_k(\boldsymbol{x}_j))$, $j=1,\ldots,N$, $k=1,\ldots,Q$, $\boldsymbol{\alpha}=(\alpha_1,\ldots,\alpha_N)^T$, $\boldsymbol{\beta}=(\beta_1,\ldots,\beta_Q)^T$, $\boldsymbol{u}=(u_1,\ldots,u_N)^T$, $\boldsymbol{0}$ is a zero vector of length $Q$, and $O$ is a $Q\times Q$ zero matrix. The system \eqref{eq:24} is obviously solvable if the coefficient matrix on the left-hand side is invertible (or non-singular) \cite{wen05}. 

\begin{definition}
The points $X=\{\boldsymbol{x}_1,\dots ,\boldsymbol{x}_N\}\subseteq {{\mathbb R}}^d$ with $N\ge Q={\dim{\pi }_m({{\mathbb R}}^d)}$ are called ${\pi }_m({{\mathbb R}}^d)$-unisolvent if the zero polynomial is the only polynomial from ${\pi }_m({{\mathbb R}}^d)$ that vanishes on all of them.
\end{definition}

\begin{theorem}\label{t11}
Suppose that $\phi$ is conditionally positive definite of order $m$ and $X$ is a ${\pi }_{m-1}({{\mathbb R}}^d)$-unisolvent set of centers. Then the system (\ref{eq:24}) is uniquely solvable.
\end{theorem}

However, the polynomial $\sum^Q_{k=1}{{\beta }_kp_k(\boldsymbol{x})}$ in \eqref{eq:21} is generally required when $\phi$ is conditionally positive definite, i.e. when $\phi$ has a polynomial growth towards infinity such as MultiQuadric and Thin Plate Spline. These functions are globally supported and conditionally positive definite in $\mathbb{R}^d$ for any $d$ of order one and $m=\nu+1$, with $\nu \in \mathbb{N}$, respectively. Therefore, in these cases the addition of a polynomial term of such orders in \eqref{eq:21} -- together with side conditions \eqref{eq:23} -- is necessary in order to guarantee the existence of a unique solution of the considered system \eqref{eq:24}. On the contrary, adding a polynomial is not usually essential when one works with positive definite RBFs such as Gaussian, Inverse MultiQuadric and Mat$\acute{\text{e}}$rn functions. This fact holds for any dimension $d$ of the space. Moreover, since all these functions are globally supported, the interpolation matrix is full and sometimes -- for some particular choices of the RBF shape parameter -- it may turn out to be very ill-conditioned \cite{fass07}.

To improve the conditioning of the system of collocation equations, CSRBFs can also be applied, even if the CSRBFs vanish beyond a user defined threshold distance $\sigma$. Therefore, only the entries in the collocation matrix corresponding to collocation nodes lying closer than $\sigma$ to a given CSRBF center are nonzero, leading to a sparse matrix. In practice, the interest in CSRBFs waned slightly as it became evident that, in order to obtain a good accuracy, the overlap distance $\sigma$ should cover most nodes in the point set, thus resulting in a populated matrix again \cite{mao08}. A typical example of family of CSRBFs is given by Wendland functions, whose support is $\left[0,\sigma\right]$, with $\sigma = 1/\epsilon$. In particular, the Wendland functions shown in Table \ref{tab_rbf} are positive definite in $\mathbb{R}^d$ for $d\leq 3$~\cite{wen95}.

In a similar representation as \eqref{eq:21}, for any linear partial differential operator$\ {\mathcal L}$, ${\mathcal L}u$ may be approximated by \cite{deh07}
\begin{equation*} 
{\mathcal L}u(\boldsymbol{x})\simeq \sum^{N}_{j=1}{{\alpha }_{j}{\mathcal L}\phi( ||\boldsymbol{x}-\boldsymbol{x}_{j}||)}+\sum^{Q}_{k=1}{{\beta }_{k}{\mathcal L}{p}_{k}(\boldsymbol{x})}.                                                   
\end{equation*} 


\subsection{RBF collocation by positive definite functions} \label{sec22}
As said in Subsection \ref{sec21}, when $\phi$ is a positive definite function, Eq. \eqref{eq:21} can be written without the additional polynomial $\ \sum^Q_{k=1}{{\beta }_kp_k(\boldsymbol{x})}$ because the solvability of the resulting interpolation system is guaranteed. So, in this case, the RBF interpolant assumes the following simplified form
\begin{equation} \label{eq:21a} 
s_{u,X}(\boldsymbol{x})=\sum_{j=1}^N \alpha_j \phi(||\boldsymbol{x}-\boldsymbol{x}_j||),
\end{equation}
where the coefficients $\alpha_{1},\ldots , \alpha_{N}$ are determined by enforcing the interpolation conditions 
\begin{equation*}
s_{u,X}(\boldsymbol{x}_{i} ) = u(\boldsymbol{x}_{i} ),\qquad i = 1,\ldots, N.
\end{equation*}
Imposing these conditions, we obtain a symmetric linear system of equations 
\begin{equation}\label{e21b}
A\boldsymbol{\alpha}=\boldsymbol{u},
\end{equation}
where $A_{ij}=\phi(||\boldsymbol{x}_{i}-\boldsymbol{x}_{j}||)$, $i,j=1,\ldots ,N$, $\boldsymbol{\alpha}= (\alpha_1, \ldots, \alpha_N)^T$, and $\boldsymbol{u} =(u_1, \ldots , u_N)^T$. If $\phi$ is a positive definite function, then the associated matrix $A$ is invertible and the RBF interpolation problem is well-posed, hence a solution to the problem exists and is unique~\cite{fass07}.

Therefore, once the vector $\boldsymbol{\alpha}$ is found, we can evaluate the RBF interpolant at any point $\boldsymbol{x}$ as
	\begin{equation}\label{e21c}
	s_{u,X}(\boldsymbol{x}) = \boldsymbol{\phi}^T(\boldsymbol{x}) \boldsymbol{\alpha},
	\end{equation}
	where $\boldsymbol{\phi}^T(\boldsymbol{x}) = (\phi (||  \boldsymbol{x} - \boldsymbol{x}_1  ||),\ldots,\phi (||  \boldsymbol{x} - \boldsymbol{x}_N  ||))$. In particular, the interpolant $s_{u,X}$ in \eqref{eq:21a} (or \eqref{e21c}) is a function of the \textsl{native Hilbert space} ${\cal N}_{\phi}(\Omega)$ uniquely associated with the RBF, and, if $u\in{\cal N}_{\phi}$, then $s_{u,X}$ is the ${\cal N}_{\phi}$-projection of $u$ into the subspace ${\cal N}_{\phi}(X)= \{\phi(||\boldsymbol{x}-\boldsymbol{x}_j||), \boldsymbol{x}_j \in X\}$~\cite{wen05}.

In the sequel, we now express the interpolant \eqref{e21c} in Lagrange form, using cardinal basis functions \cite{saf15}. The cardinal basis functions $\psi_{j}(\boldsymbol{x})$, $j = 1,\ldots, N$, have the property
\begin{equation*}
\psi_{j} (\boldsymbol{x}_{i})= \left\lbrace \begin{array}{ll}
1 & \text{if} ~i=j,\\
0 & \text{if} ~i\neq j,
\end{array}  \right. j=1,\ldots,N,
\end{equation*}
leading to the alternative formulation for the interpolant
\begin{equation}\label{e21d}
s_{u,X}(\boldsymbol{x})=\boldsymbol{\psi}^T(\boldsymbol{x})\boldsymbol{u}
\end{equation}
where $\boldsymbol{\psi}^T(\boldsymbol{x})=(\psi_{1}(\boldsymbol{x}), \ldots, \psi_{N}(\boldsymbol{x}))$. Combining (\ref{e21c}), (\ref{e21d}), and (\ref{e21b}), we obtain 
\begin{equation}\label{e21e}
s_{u,X}(\boldsymbol{x})=\boldsymbol{\psi}^T(\boldsymbol{x})\boldsymbol{u}=\boldsymbol{\phi}^T(\boldsymbol{x})\boldsymbol{\alpha}=\boldsymbol{\phi}^T(\boldsymbol{x})A^{-1}\boldsymbol{u},
\end{equation}
so that we may deduce the following relation between the cardinal basis and the original radial basis
\begin{equation*}
\boldsymbol{\psi}^T(\boldsymbol{x})=\boldsymbol{\phi}^T(\boldsymbol{x})A^{-1}.
\end{equation*}
This transformation is valid whenever the matrix $A$ is nonsingular, i.e. for distinct data points $\boldsymbol{x}_{1},\ldots, \boldsymbol{x}_{N} \in X$ and positive definite RBFs.

For a linear operator $\mathcal{L}$, we have
\begin{equation*}
\mathcal{L}s_{u,X}(\boldsymbol{x})=\sum_{j=1}^{N}\mathcal{L}\psi_{j}(\boldsymbol{x})u(\boldsymbol{x}_{j}).
\end{equation*}
To evaluate $\mathcal{L}s_{u,X}(\boldsymbol{x})$ at the nodes, i.e., to evaluate $\boldsymbol{s}_{\mathcal{L}} =(\mathcal{L}s_{u,X}(\boldsymbol{x}_{1}),\ldots,\mathcal{L}s_{u,X}(\boldsymbol{x}_{N}))^{T}$, we need the differentiation matrix $\Psi_{\mathcal{L}}^T =(\mathcal{L}\psi_{j} (\boldsymbol{x}_{i}))$, $i,j=1,\ldots,N$. From relation (\ref{e21e}), we acquire
\begin{equation*}
\boldsymbol{s}_{\mathcal{L}}=\Psi_{\mathcal{L}}^T \boldsymbol{u}=\Phi_{\mathcal{L}}^T A^{-1}\boldsymbol{u},
\end{equation*}
where $\Phi_{\mathcal{L}}^T = (\mathcal{L}\phi (||\boldsymbol{x}_i-\boldsymbol{x}_{j}||))$, $i,j=1,\ldots,N$.

Therefore, in order to solve a time-dependent PDE problem, we can use an RBF interpolant expressed in the Lagrangian form, approximating the solution $u(\boldsymbol{x},t)$ as follows
\begin{align}\label{lagtdi}
s_{u,X}(\boldsymbol{x},t) =\sum_{j=1}^{N} \psi_j(\boldsymbol{x}) u_j(t), \qquad t\geq 0,
\end{align}
where $u_j(t)\approx u(\boldsymbol{x}_j,t)$ are the unknown functions to be determined.


\section{PUM collocation by RBFs}\label{sec3}

In this section, we consider the RBF-PUM collocation method for PDEs, expressing this local meshfree method in terms of its weight functions and local RBF interpolants. In the following, for simplicity we refer to the case of PUM collocation by positive definite RBFs, but the whole presentation could be done for conditionally positive definite RBFs, as outlined in Subsection \ref{sec21}.


\subsection{RBF-PUM for interpolation problems}
Let $\Omega \subseteq \mathbb{R}^d$ be an open bound\-ed domain, and let $\{\Omega_j\}_{j=1}^{M}$ be an open and bounded covering of $\Omega$ satisfying some mild overlap condition among the subdomains (or patches) $\Omega_j$, i.e. the overlap among the subdomains must be sufficient so that each interior point $x \in \Omega$ is located in the interior of at least one subdomain $\Omega_{j}$. Moreover, the set $I(\boldsymbol{x}) = \{j : \boldsymbol{x} \in \Omega_j \}$ is uniformly bounded by the constant $K$ -- independent of $M$ -- on $\Omega$, namely ${card}(I(\boldsymbol{x}))\leq K$, for $\boldsymbol{x} \in \Omega$, where $ \Omega  \subseteq \bigcup_{j=1}^{M} \Omega_j$. Associated with the subdomains we define a partition of unity $\{w_j\}_{J=1}^{M}$ subordinated to the covering $\{\Omega_j\}_{j=1}^{M}$ such that 
$$ 
\sum_{j = 1}^M w_j( \boldsymbol{x})=1, \qquad \boldsymbol{x} \in \Omega,
$$
where the weight function $w_j: \Omega_j \rightarrow \mathbb{R}$ is a compactly supported, nonnegative and continuous with ${supp}(w_j)  \subseteq \Omega_j$. For each subdomain we may thus construct, similarly to~\eqref{eq:21a}, a local RBF interpolant $s_{u_j,X_j}:\Omega_j \rightarrow \mathbb{R}$ of the form 
\begin{align} \label{loc_rbf}
	s_{u_j,X_j}( \boldsymbol{x}) = \sum_{i=1}^{N_j} \alpha_i^j \phi (||  \boldsymbol{x} - \boldsymbol{x}_i^j  ||_2),
\end{align}
where $N_j$ is the number of nodes in $\Omega_j$, i.e. the points $\boldsymbol{x}_i^j  \in X_{j} = X \cap \Omega_j$. Therefore the global RBF-PUM interpolant is defined as
	\begin{equation} \label{e1}
	{\cal P}_{u,X}( \boldsymbol{x})= \sum_{j=1}^{M} w_j ( \boldsymbol{x}) s_{u_j,X_j}( \boldsymbol{x} ), \qquad \boldsymbol{x} \in \Omega.
	\end{equation}
If the functions $s_{u_j,X_j}$, $j=1,\ldots,M$, satisfy the interpolation conditions 
\begin{align} \label{loc_rbf_int}
s_{u_j,X_j}( \boldsymbol{x}_i^j )= u( \boldsymbol{x}_i^j), \qquad \boldsymbol{x}_i^j \in \Omega_j,\qquad i=1,\ldots,N_j,
\end{align} 
then the global interpolant \eqref{e1} inherits the interpolation property of the local interpolants, i.e.
		\begin{align} \label{iip}
		{\cal P}_{u,X}( \boldsymbol{x}_i^j ) = \sum_{j=1}^{M} w_j ( \boldsymbol{x}_i^j) s_{u_j,X_j}( \boldsymbol{x}_i^j ) = \sum_{j=1}^M u( \boldsymbol{x}_i^j ) w_j ( \boldsymbol{x}_i^j) = u( \boldsymbol{x}_i^j).
		\end{align}
Solving the $j$-th interpolation problem \eqref{loc_rbf_int} leads to the local RBF linear system
\begin{align*} 
 A_{j}\boldsymbol{\alpha}_{j}=\boldsymbol{u}_{j},
\end{align*}
where $A_j$ ia a ${N_j\times N_j}$ matrix of entries $A_{ik}^j=\phi(||\boldsymbol{x}_i^j-\boldsymbol{x}_k^j||)$, $i,k=1,\ldots,N_j$, $\boldsymbol{\alpha}_j=(\alpha_1^j,\ldots,\alpha_{N_j}^j)^T$, and $\boldsymbol{u}_j=(u_1^j,\ldots,u_{N_j}^j)^T$. Note that -- similarly to \eqref{e21b} -- existence and uniqueness of the solution and nonsingularity of the local matrix $A_{j}$ is guaranteed by the use of positive definite RBFs $\phi$~\cite{fass07}.

The PUM weight function $w_j$ can be constructed using the Shepard method given by
	\begin{align} \label{ShepWei}
	w_j(\boldsymbol{x}) = \frac{\varphi_j(\boldsymbol{x})}{\sum_{k \in I(\boldsymbol{x})} \varphi_k(\boldsymbol{x})}, \qquad j=1,\ldots,M,
	\end{align}
	$\varphi_j(\boldsymbol{x})$ being a compactly supported function with support on $\Omega_j$ such as the Wendland $C^2$ function shown in Table~\ref{tab_rbf}. Such functions are scaled with a shape parameter $\gamma$ to get $\varphi_j(\boldsymbol{x})=\varphi(\gamma||\boldsymbol{x} - \boldsymbol{\xi}_j||)$, where $\boldsymbol{\xi}_j$ is the center of the weight function. It follows that $w_{j} (\boldsymbol{x}) = 0,~ \forall j \notin I(\boldsymbol{x})$. Therefore, we may rewrite Eq. \eqref{e1} as
\begin{equation}\label{e4}
\mathcal{P}_{u,X}(\boldsymbol{x})=\sum_{j \in I(\boldsymbol{x})}w_{j}(\boldsymbol{x})s_{u_j,X_j}(\boldsymbol{x}).
\end{equation}
and the weight functions $w_j(\boldsymbol{x})$ in \eqref{ShepWei} satisfy the partition of unity property
$$ 
\sum_{j \in I(\boldsymbol{x})} w_j(\boldsymbol{x})=1.
$$
Hence, equivalently to Eqs. \eqref{loc_rbf_int} and \eqref{iip}, if the local fits $s_{u_j,X_j}$ in Eq. \eqref{e4} interpolate at a given data point $\boldsymbol{x}_i$, i.e. $s_{u_j,X_j}(\boldsymbol{x}_{i}) = u(\boldsymbol{x}_{i})$ for each node $\boldsymbol{x}_{i} \in \Omega_{j}$, then the global fit also interpolates at that point, that is 
\begin{equation*}
\mathcal{P}_{u,X}(\boldsymbol{x}_{i})=\sum_{j \in I(\boldsymbol{x}_{i})}w_{j}(\boldsymbol{x}_{i})s_{u_j,X_j}(\boldsymbol{x}_{i})=u(\boldsymbol{x}_{i})\sum_{j \in I(\boldsymbol{x}_{i})}w_{j}(\boldsymbol{x}_{i})=u(\boldsymbol{x}_{i}).
\end{equation*}

To be able to formulate error bounds, we first consider some technical conditions and then define a few assumptions on regularity of $\Omega_j$~\cite{wen02}. So we require the partition of unity functions $w_j$ to be \textsl{k-stable}. This means that each $w_j \in C^k(\mathbb{R}^d)$ satisfies, for every multi-index $\boldsymbol{\mu} \in \mathbb{N}_0^d$ with $|\boldsymbol{\mu}| \leq k$, the inequality
\begin{equation}
	\left\|D^{\boldsymbol{\mu}}w_j\right\|_{L_{\infty}(\Omega_j)}\leq C_{\boldsymbol{\mu}}/\delta_j^{|\boldsymbol{\mu}|}. \nonumber
\end{equation}
where $C_{\boldsymbol{\mu}}$ is some positive constant, and $\delta_j$ = diam($\Omega_j$).

	\begin{definition} \label{defpr}
	Let $\Omega \subseteq  \mathbb{R}^d$ be bounded and let $\boldsymbol{x}_1,\ldots,\boldsymbol{x}_N \in X \subseteq \Omega$ be given. An open and bounded covering $\{\Omega_j\}_{j=1}^{M}$ is called regular for $(\Omega,X)$ if the following conditions hold:
	\begin{itemize}
		\item
		for each $\boldsymbol{x} \in \Omega$, the number of subdomains $\Omega_j$ with $\boldsymbol{x} \in \Omega_j$ is bounded by a global constant $K$;
		\item
		each subdomain $\Omega_j$ satisfies an interior cone condition~\cite{wen05};
		\item
		the local fill distances $h_{X_j, \Omega_j}$ are uniformly bounded by the global fill distance
		\begin{align*}
	h_{ X, \Omega} =  \sup_{ \boldsymbol{x} \in \Omega} \min_{ \boldsymbol{x}_i  \in X} || \boldsymbol{x} - \boldsymbol{x}_i||.
	\end{align*}
	\end{itemize}
\end{definition}

	After defining the space $C_{ \nu}^{k}  ( \mathbb{R}^d ) $ of all functions $u \in C^k$ whose derivatives of order $ | \boldsymbol{\mu} |=k $ satisfy $ D^{ \boldsymbol{\mu}} u ( \boldsymbol{x} ) = {\cal O} ( || \boldsymbol{x} ||^{ \nu} ) $ for $ || \boldsymbol{x} ||_2 \rightarrow 0$, we consider the following convergence result \cite{wen05}.

\begin{theorem} \label{theores}
	Let $\Omega \subseteq  \mathbb{R}^d$ be open and bounded and $\boldsymbol{x}_1,\ldots,\boldsymbol{x}_N \in X\subseteq \Omega$. Let $\phi \in C_{\nu}^k(\mathbb{R}^d)$ be a conditionally positive definite function of order $m$. If $\{\Omega_j\}_{j=1}^{M}$ is a regular covering for $(\Omega, X)$ and $\{w_j\}_{j=1}^{M}$ is $k$-stable for $\{\Omega_j\}_{j=1}^{M}$, then the error between $u \in {\cal N}_{\phi}(\Omega)$ and its PUM interpolant ${\cal P}_{u,X}$ is bounded by
	\begin{equation}
	|D^{\boldsymbol{\mu}}u(\boldsymbol{x}) - D^{\boldsymbol{\mu}}{\cal P}_{u,X}(\boldsymbol{x})| \leq C h_{X, \Omega}^{(k+\nu)/2 - |\boldsymbol{\mu}|} |u|_{{\cal N}_{\phi}(\Omega)}, \nonumber
	\end{equation}
	for all $\boldsymbol{x} \in \Omega$, and all $|\boldsymbol{\mu}| \leq k/2$.
\end{theorem}

If we compare the convergence result reported in Theorem \ref{theores} with the global error estimates in~\cite{wen05}, we can note that the PUM interpolant preserves the local approximation order for the global fit. So we can efficiently compute large RBF interpolants problem by solving many small RBF interpolation problems and then glue them together with the global partition of unity. It follows that the PUM-based approach is a simple and computationally efficient tool to decompose a large interpolation problem into many small problems, while at the same time ensuring that the accuracy obtained for the local fits is carried over to the global one.   


\subsection{RBF-PUM for differential problems}
When a time-dependent PDE is numerically solved by a RBF-PUM scheme, the solution $u(\boldsymbol{x},t)$ of the differential problem is approximated by the global interpolant
\begin{equation}\label{e5}
\mathcal{P}_{u,X}(\boldsymbol{x},t)=\sum_{j \in I(\boldsymbol{x})}w_{j}(\boldsymbol{x})s_{u_j,X_j}(\boldsymbol{x},t), \qquad t\geq 0,
\end{equation}
where -- similarly to Eq. \eqref{lagtdi} -- $s_{u_j,X_j}(\boldsymbol{x},t)$ is a local RBF interpolant defined on $\Omega_{j}$ of the form
\begin{equation}\label{e6}
s_{u_j,X_j}(\boldsymbol{x},t)=\sum_{k\in J(\Omega_{j})} \psi_{k}(\boldsymbol{x})u_{k}(t),
\end{equation}
with $J(\Omega_{j}) = \lbrace k : \boldsymbol{x}_{k} \in \Omega_{j} \rbrace$ that defines the set of nodes in $\Omega_{j}$. From Eqs. \eqref{e5} and \eqref{e6}, we obtain that the RBF-PUM interpolant may be expressed as follows
\begin{equation}\label{e7}
\mathcal{P}_{u,X}(\boldsymbol{x},t)=\sum_{j \in I(\boldsymbol{x})}w_{j}(\boldsymbol{x})\sum_{k\in J(\Omega_{j})} \psi_{k}(\boldsymbol{x})u_{k}(t)=\sum_{j \in I(\boldsymbol{x})} \sum_{k\in J(\Omega_{j})} \left( w_{j}(\boldsymbol{x}) \psi_{k}(\boldsymbol{x}) \right)u_{k} (t).
\end{equation}
Hence, if we interpolate the initial condition of a time-dependent problem, we have $\mathcal{P}_{u,X}(\boldsymbol{x}_k,0) = u(\boldsymbol{x}_k,0)$ for all $k$, while $\mathcal{P}_{u,X}(\boldsymbol{x}_k,t) \approx u(\boldsymbol{x}_k,t)$ for $t>0$ \cite{saf15}.

In a PDE problem using the RBF-PUM approximation scheme \eqref{e7} requires to compute the effect of applying a spatial differential operator $\mathcal{L}$ at the interior data points. Denoting first by $\boldsymbol{\mu}$ and $\boldsymbol{\nu}$ the multi-indices for usual rules and using then Leibniz’s rule, we can compute a derivative term of order $\boldsymbol{\mu}$ of the global fit \eqref{e7}, obtaining the derivative rule
\begin{equation}\label{e8}
\begin{split}
\dfrac{\partial^{\vert \boldsymbol{\mu} \vert}}{\partial \boldsymbol{x}^{\boldsymbol{\mu}}} \mathcal{P}_{u,X}(\boldsymbol{x},t) & = \sum_{j \in I(\boldsymbol{x})} \sum_{k\in J(\Omega_{j})}  \dfrac{\partial^{\vert \boldsymbol{\mu} \vert}}{\partial \boldsymbol{x}^{\boldsymbol{\mu}}} \left( w_{j}(\boldsymbol{x}) \psi_{k}(\boldsymbol{x}) \right) u_{k}(t)\\
& = \sum_{j \in I(\boldsymbol{x})} \sum_{k\in J(\Omega_{j})}  \left( \sum_{\boldsymbol{\nu} \leq \boldsymbol{\mu}} \left(\begin{array}{l}\boldsymbol{\mu}\\\boldsymbol{\nu}\end{array}\right) \dfrac{\partial^{\vert \boldsymbol{\mu} - \boldsymbol{\nu} \vert} w_{j}}{\partial \boldsymbol{x}^{\boldsymbol{\mu} - \boldsymbol{\nu}}} (\boldsymbol{x}) \dfrac{\partial^{\vert \boldsymbol{\nu} \vert} \psi_{k}}{\partial \boldsymbol{x}^{\boldsymbol{\nu}}} (\boldsymbol{x}) \right)u_{k}(t).
\end{split}
\end{equation}
If we fix $\boldsymbol{x} = \boldsymbol{x}_{i}$ and $k$ in Eq. (\ref{e8}), we get the $ik$-element of the global differentiation matrix corresponding to the $\boldsymbol{\mu}$-derivative. Moreover, when having composite linear operators, we take the contributions from each term.


\section{RBF-PUM collocation method based on finite difference} \label{cd_pp_prbs}
In this section, we apply RBF partition of unity collocation method based on finite difference for two time-dependent PDEs. Specifically, we present our collocation scheme to firstly solve an unsteady convection-diffusion equation in Subsection \ref{cd-p}, and then a pseudo-parabolic equation in Subsection \ref{pp-p}.


\subsection{Collocation scheme for a convection-diffusion problem} \label{cd-p}
Let us consider the following unsteady convection-diffusion equation
\begin{equation}\label{e9}
\dfrac{\partial u(\boldsymbol{x},t)}{\partial t}=\kappa \vartriangle u(\boldsymbol{x},t) + \nu \cdot \triangledown u(\boldsymbol{x},t) \equiv \mathcal{L} u(\boldsymbol{x},t),\qquad \boldsymbol{x} \in \Omega \subset \mathbb{R}^{d},\quad t>0,
\end{equation}
where $\mathcal{L}$ is the convection-diffusion operator, $\vartriangle$ and $\triangledown$ denote the Laplacian and the gradient operator, respectively, $\kappa$ is the diffusion coefficient, $\nu$ is a constant velocity vector, and $u(\boldsymbol{x},t)$ may represent concentration or temperature for mass or heat transfer. The Eq. \eqref{e9} must be supplemented with an initial condition of the form
\begin{equation}\label{e17}
u(\boldsymbol{x},0)=u_{0}(\boldsymbol{x}),
\end{equation}
and boundary conditions
\begin{equation}\label{e14}
\mathcal{B}u(\boldsymbol{x},t)=g(\boldsymbol{x},t), \qquad \boldsymbol{x} \in \partial \Omega, \quad t>0,
\end{equation}
where $g(\boldsymbol{x},t)$ is a known function, and $\mathcal{B}$ is a boundary operator, which can be of Dirichlet, Neumann or mixed type, and $\partial \Omega$ denotes the boundary of $\Omega$. In the case of Dirichlet boundary conditions, we may discretize time derivative of the PDE \eqref{e9} by usual finite difference formula, using the following $\theta$-weighted scheme
\begin{equation}\label{e10}
\dfrac{u^{n+1}(\boldsymbol{x})-u^{n}(\boldsymbol{x})}{\delta t}= \theta \left( \kappa \vartriangle u^{n+1}(\boldsymbol{x})+\nu \cdot \triangledown u^{n+1}(\boldsymbol{x}) \right) + (1-\theta) \left( \kappa \vartriangle u^{n}(\boldsymbol{x})+\nu \cdot \triangledown u^{n}(\boldsymbol{x})  \right),
\end{equation}
or equivalently
\begin{equation}\label{e10b}
\dfrac{u^{n+1}(\boldsymbol{x})-u^{n}(\boldsymbol{x})}{\delta t}= \theta \mathcal{L} u^{n+1}(\boldsymbol{x}) + (1-\theta) \mathcal{L} u^{n}(\boldsymbol{x}),
\end{equation}
where $0\leq \theta \leq 1$, $u^{n+1}(\boldsymbol{x})=u(\boldsymbol{x},t^{n+1})$, $t^{n+1}=t^{n}+\delta t$, and $\delta t$ is the time step size. Rearranging Eq. \eqref{e10} (or Eq. \eqref{e10b}), we obtain
\begin{equation}\label{e13}
u^{n+1}(\boldsymbol{x})+\eta \left( \kappa \vartriangle u^{n+1}(\boldsymbol{x})+\nu \cdot \triangledown u^{n+1}(\boldsymbol{x})  \right)=u^{n}(\boldsymbol{x})+\zeta \left( \kappa \vartriangle u^{n}(\boldsymbol{x})+\nu \cdot \triangledown u^{n}(\boldsymbol{x})  \right),
\end{equation}
or
\begin{equation}\label{e13b}
u^{n+1}(\boldsymbol{x})+\eta \mathcal{L} u^{n+1}(\boldsymbol{x}) =u^{n}(\boldsymbol{x})+\zeta \mathcal{L} u^{n}(\boldsymbol{x}),
\end{equation}
where $\eta=- \theta \delta t $ and $\zeta=(1-\theta) \delta t$.

Referring now to Eq. \eqref{e7}, we get that $u^{n}(\boldsymbol{x})$ can be approximated by
\begin{equation}\label{e11}
u^{n}(\boldsymbol{x})\simeq p^{n}(\boldsymbol{x})= \sum_{j \in I(\boldsymbol{x})} \sum_{k\in J(\Omega_{j})}  \left( w_{j}(\boldsymbol{x}) \psi_{k}(\boldsymbol{x}) \right)u_{k}^{n}.
\end{equation}
If we assume that $\mathfrak{I}$ and $\mathfrak{B}$ denote the indexes of internal and boundary points, respectively, and suppose that $N$ is the total number of centers, i.e. $N=N_{\mathfrak{I}} +N_{\mathfrak{B}}$, then the $N\times N$ matrix $A$ can be split into two matrices $A_{\mathfrak{I}}$ and $A_{\mathfrak{B}}$ so that
\begin{equation} \label{mat-A}
A=A_{\mathfrak{I}}+A_{\mathfrak{B}}, 
\end{equation}
where
\begin{equation*}
\begin{split}
&A= \left[ w_{j}(\boldsymbol{x}_{i})\psi_{k}(\boldsymbol{x}_{i})~ \text{for}~ (j \in I(\boldsymbol{x}_{i}),k\in J(\Omega_{j}),i=1,\ldots,N)~ \text{and}~ 0 ~\text{elsewhere} \right]_{N\times N},\\
&A_{\mathfrak{I}} = \left[a_{ij}~\text{for}~(i\in \mathfrak{I}, 1\leq j\leq N)~\text{and}~0~\text{elsewhere} \right],\\
&A_{\mathfrak{B}} = \left[a_{ij}~\text{for}~(i\in \mathfrak{B}, 1\leq j\leq N)~\text{and}~0~\text{elsewhere} \right].
\end{split}
\end{equation*}
Using the notation $\mathcal{L}A$ to designate the matrix of the same dimension as $A$, whose elements are $\widehat{a}_{ij}=\mathcal{L}a_{ij},~1\leq i,j\leq N$, and substituting (\ref{e11}) in Eq. \eqref{e13} (or Eq. \eqref{e13b}) together with (\ref{e14}), we can write the resulting sparse system, in the matrix form, as
\begin{equation}\label{e15}
C \boldsymbol{u}^{n+1}= D \boldsymbol{u}^{n}+\boldsymbol{v}^{n+1},
\end{equation}
where
\begin{equation*}
\begin{split}
 &C= A+\eta \kappa \vartriangle A_{\mathfrak{I}} + \eta \nu\cdot \triangledown A_{\mathfrak{I}} ,\\
 &D= A_{\mathfrak{I}} + \zeta \kappa \vartriangle A_{\mathfrak{I}} + \zeta \nu\cdot \triangledown A_{\mathfrak{I}} ,\\
 &\boldsymbol{v}^{n+1}= \left[ g_{i}^{n+1}~\text{for}~(i\in \mathfrak{B})~\text{and}~0~\text{elsewhere} \right]^{T},\\
 &\boldsymbol{u}^{n}=\left( u_{1}^{n},~\ldots,~u_{N}^{n} \right)^{T}.
\end{split}
\end{equation*}
The system \eqref{e15} is obtained by combining Eq. \eqref{e13} (or Eq. \eqref{e13b}), which applies to the internal points, and Eq. \eqref{e14} that refers to the boundary points. Therefore, using the initial condition represented by Eq. \eqref{e17}, we can compute $\boldsymbol{u}^{n+1}$ by solving Eq. (\ref{e15}). Then, by substituting such values of $\boldsymbol{u}^{n}$ in 
\begin{equation} \label{wn-step}
\boldsymbol{p}^{n} = A \boldsymbol{u}^{n},
\end{equation}
we can obtain the approximated solution of the PDE at time level $n$. Nevertheless, as well as for the RBF collocation scheme \cite{hon01}, also in the RBF-PUM-FD case we do not have a theoretical proof that shows the matrix $C$ is invertible when $\theta > 0$. So, here, we deal with the problem numerically, postponing the treatment of more theoretical issues such as well-posedness to future works. For the explicit scheme, i.e. for $\theta = 0$, we only need to invert the matrix $A$, whose invertibility is guaranteed provided that the set of collocation points is distinct.


\subsection{Collocation scheme for a pseudo-parabolic problem} \label{pp-p}
We consider the following third order equation
\begin{equation}\label{e9-1}
\dfrac{\partial u(\boldsymbol{x},t)}{\partial t} - \alpha \vartriangle u(\boldsymbol{x},t) - \beta \dfrac{\partial \vartriangle u(\boldsymbol{x},t)}{\partial t}  = f(\boldsymbol{x},t),\qquad \boldsymbol{x} \in \Omega \subset \mathbb{R}^{d},\quad t>0,
\end{equation}
where $\vartriangle$ denotes the Laplacian operator, $\alpha, \beta>0$ are known constants, and $f(\boldsymbol{x},t)$ is given continuous function. The above equation is usually called Sobolev-type or pseudo-parabolic equation. To Eq. \eqref{e9-1} we add initial condition
\begin{equation}\label{e17-1}
u(\boldsymbol{x},0)=u_{0}(\boldsymbol{x}),
\end{equation}
and boundary conditions
\begin{equation}\label{e14-1}
\mathcal{B}u(\boldsymbol{x},t)=g(\boldsymbol{x},t), \qquad \boldsymbol{x} \in \partial \Omega, \quad t>0,
\end{equation}
where $g(\boldsymbol{x},t)$ is a known continuous function, and $\mathcal{B}$ is a boundary operator, which can be of Dirichlet, Neumann or mixed type, and $\partial \Omega$ denotes the boundary of $\Omega$. In the case of Dirichlet boundary conditions, we can discretize time derivative of Eq. \eqref{e9-1} via $\theta$-weighted finite difference scheme as
\begin{equation}\label{e10-1}
\dfrac{u^{n+1}(\boldsymbol{x})-u^{n}(\boldsymbol{x})}{\delta t}-\beta \dfrac{\vartriangle u^{n+1}(\boldsymbol{x})- \vartriangle u^{n}(\boldsymbol{x})}{\delta t}= \theta \left( \alpha \vartriangle u^{n+1}(\boldsymbol{x})+f^{n+1}(\boldsymbol{x}) \right) + (1-\theta) \left( \alpha \vartriangle u^{n}(\boldsymbol{x})+f^{n}(\boldsymbol{x})  \right),
\end{equation}
where $0\leq \theta \leq 1$, $u^{n+1}(\boldsymbol{x})=u(\boldsymbol{x},t^{n+1})$, $t^{n+1}=t^{n}+\delta t$, and $\delta t$ is the time step size. Eq. \eqref{e10-1} can then be rewritten in the following form
\begin{equation}\label{e13-1}
u^{n+1}(\boldsymbol{x})-(\eta+\beta)\vartriangle u^{n+1}(\boldsymbol{x})=u^{n}(\boldsymbol{x})+(\zeta-\beta)\vartriangle u^{n+1}(\boldsymbol{x})+z^{n+1}(\boldsymbol{x}),
\end{equation}
where $\eta= \theta \delta t \alpha $, $\zeta=(1-\theta) \delta t \alpha$, and $z^{n+1}=\delta t (\theta f^{n+1}(\boldsymbol{x}) +(1-\theta) f^{n}(\boldsymbol{x}) )$.

Similarly to the previous subsection, considering Eq. \eqref{e7}, $u^{n}(\boldsymbol{x})$ can be approximated as in \eqref{e11}. Therefore, always denoting by $\mathfrak{I}$ and $\mathfrak{B}$ the indexes of internal and boundary points, respectively, and $N$ being the total number of nodes so that $N=N_{\mathfrak{I}} +N_{\mathfrak{B}}$, the matrix $A \in \mathbb{R}^{N\times N}$ is given by the sum of matrices $A_{\mathfrak{I}}$ and $A_{\mathfrak{B}}$ as in \eqref{mat-A},
where
\begin{equation*}
\begin{split}
&A= \left[ w_{j}(\boldsymbol{x}_{i})\psi_{k}(\boldsymbol{x}_{i})~ \text{for}~ (j \in I(\boldsymbol{x}_{i}),k\in J(\Omega_{j}),i=1,\ldots,N)~ \text{and}~ 0 ~\text{elsewhere} \right]_{N\times N},\\
&A_{\mathfrak{I}} = \left[a_{ij}~\text{for}~(i\in \mathfrak{I}, 1\leq j\leq N)~\text{and}~0~\text{elsewhere} \right],\\
&A_{\mathfrak{B}} = \left[a_{ij}~\text{for}~(i\in \mathfrak{B}, 1\leq j\leq N)~\text{and}~0~\text{elsewhere} \right].
\end{split}
\end{equation*}
Now, substituting \eqref{e11} in Eq. \eqref{e13-1} together with \eqref{e14-1}, we obtain a sparse matrix system of the form \eqref{e15}, with
\begin{equation*}
\begin{split}
 &C= A-(\eta+\beta)\vartriangle A_{\mathfrak{I}},\\
 &D= A_{\mathfrak{I}}+(\zeta-\beta)\vartriangle A_{\mathfrak{I}},\\
 &\boldsymbol{v}^{n+1}= \left[ z_{i}^{n+1}~\text{for}~(i\in \mathfrak{I})~\text{and}~g_{j}^{n+1}~\text{for}~(j\in \mathfrak{B}) \right]^{T},\\
 &\boldsymbol{u}^{n}=\left( u_{1}^{n},~\ldots,~u_{N}^{n} \right)^{T}.
\end{split}
\end{equation*}
The associated system \eqref{e15} is obtained by combining Eq. \eqref{e13-1}, which refers to the internal points, and Eq. \eqref{e14-1} that applies to the boundary points. Therefore, using the initial condition represented by Eq. \eqref{e17-1}, we solve Eq. (\ref{e15}) computing $\boldsymbol{u}^{n+1}$. Finally, replacing such values of $\boldsymbol{u}^{n}$ in \eqref{wn-step}, we get solution of the PDE at time level $n$. In conclusion, however, all the remarks about the invertibility issue of the matrix $C$ done previously for the convection-diffusion equation can be extended to this differential problem as well. 


\section{Stability analysis}\label{sec5}
In this section, as in \cite{chi06}, we provide a numerical stability analysis of the RBF-PUM-FD collocation scheme considering explicitly the two time-dependent equations studied. In so doing, we introduce a perturbation in Eq. \eqref{e15}, i.e. we set $\boldsymbol{e}^{n}=\boldsymbol{u}^{n}-\boldsymbol{p}^{n}$, where $\boldsymbol{u}^{n}$ is the discrete exact solution and $\boldsymbol{p}^{n}$ is the approximated numerical solution. We can then express the equation for the error $\boldsymbol{e}^{n+1}$ as
\begin{equation}\label{e18}
\boldsymbol{e}^{n+1}=K\boldsymbol{e}^{n},
\end{equation}
where $K=A{C}^{-1}{D}A^{-1}$ is the so-called amplification matrix. It follows that the numerical scheme is stable if the error $\boldsymbol{e}^{n} \rightarrow 0$, as $n\rightarrow \infty$. This condition of numerical stability is guaranteed provided that the spectral radius $\rho$ of the matrix $K$ is such that $\rho ({K})\leq 1$. Substituting ${K}$ in Eq. \eqref{e18} we obtain
\begin{equation}\label{e19}
{C}A^{-1} e^{n+1}={D}A^{-1} e^{n}.
\end{equation}
Now, assuming Dirichlet boundary conditions, Eq. (\ref{e19}) can be expressed in the form
\begin{equation}\label{e20}
\underbrace{\left[ {I}-\theta \delta t {M} \right]}_{P} \boldsymbol{e}^{n+1}=\underbrace{\left[ {I}+(1-\theta) \delta t {M} \right]}_{Q} \boldsymbol{e}^{n},
\end{equation}
where ${I} \in \mathbb{R}^{N\times N}$ is the identity matrix and the matrix ${M}=\mathcal{L} A_{\mathfrak{I}} A^{-1}$.

So referring to Eq. \eqref{e20}, we can deduce that stability is assured if all the eigenvalues of the matrix $P^{-1} Q$ are less than one. This fact occurs when the inequality
\begin{equation}\label{e21}
\left\vert \dfrac{1+(1-\theta) \delta t \lambda_{M}}{1-\theta \delta t \lambda_{M}} \right\vert \leq 1,
\end{equation}
$\lambda_{M}$ being an eigenvalue of matrix ${M}$. Such eigenvalues can be calculated by solving the following generalized eigenvalue problem
\begin{equation*}
\mathcal{L} A_{\mathfrak{I}} \boldsymbol{s}= \lambda_{M} A \boldsymbol{s}.
\end{equation*}
In case of Crank-Nicholson scheme, i.e. taking $\theta=1/2$, the inequality (\ref{e21}) is always satisfied if $\lambda_{M} \leq 0$, as well as for $\theta = 1$. This means that, when $\lambda_{M} \leq 0$, the numerical scheme is unconditionally stable. If, instead, we consider the explicit scheme for $\theta=0$, the stability condition becomes
\begin{equation*}
\left\vert 1 + \delta t \lambda_{M}  \right\vert \leq 1,
\end{equation*}
so the associated method turns out to be stable if
\begin{equation*}
\delta t \leq -\dfrac{2}{\lambda_{M}}~\text{and}~\lambda_{M} \leq 0.
\end{equation*} 

In case of Eq. \eqref{e9-1} the stability analysis is rather similar to that outlined eaelier. Adding indeed a perturbation, the error assumes the form \eqref{e18}, where $K$ can be split as to obtain Eq. \eqref{e19}. However, in this case Eq. \eqref{e19} assumes the form
\begin{equation}\label{e20-1}
\underbrace{\left[ {I}-(\eta+\beta) {M} \right]}_{{P}} \boldsymbol{e}^{n+1}=\underbrace{\left[ {I}+(\zeta-\beta) {M} \right]}_{{Q}} \boldsymbol{e}^{n},
\end{equation}
where ${I}$ is the identity matrix, ${M}=\vartriangle A_{\mathfrak{I}} A^{-1}$, $\eta= \theta \delta t \alpha $ and $\zeta=(1-\theta) \delta t \alpha$. From \eqref{e20-1} it follows that 
\begin{equation*}
\boldsymbol{e}^{n+1}=P^{-1}Q\boldsymbol{e}^{n},
\end{equation*}
and therefore the numerical scheme is stable provided that $\rho(P^{-1}Q) \leq 1$. This condition is thus true when
\begin{equation}\label{e21-1}
\left\vert \dfrac{1+\left((1-\theta) \delta t \alpha - \beta\right) \lambda_{M}}{1-\left(\theta \delta t \alpha + \beta\right) \lambda_{M}} \right\vert \leq 1,
\end{equation}
where $\lambda_{M}$ denotes an eigenvalue of ${M}$. In particular, we observe again that for $\theta = 1/2$ and $\theta = 1$ the inequality \eqref{e21-1} is always satisfied if $\lambda_M \leq 0$.

As shown in \eqref{e21} and \eqref{e21-1}, the stability of our method depends at least on three factors, i.e. $\theta$, $\delta t$ and $\lambda_{M}$. Furthermore, from Table \ref{tab_rbf} we notice that RBFs are usually dependent from a shape parameter $\epsilon$, which may influence the numerical stability of the collocation schemes. In Section \ref{sec6} we will investigate numerical stability analysis for the two benchmark differential problems discussed in this paper.


\section{Numerical experiments}\label{sec6}

In this section we report the performance of our RBF-PUM-FD method which is measured through numerical experiments. All these results illustrated in some tables and figures have been carried out in \textsc{Matlab} on a laptop with a 2.6 GHz Intel Core i5 processor.

\begin{figure}
\centering
\includegraphics[scale=0.6]{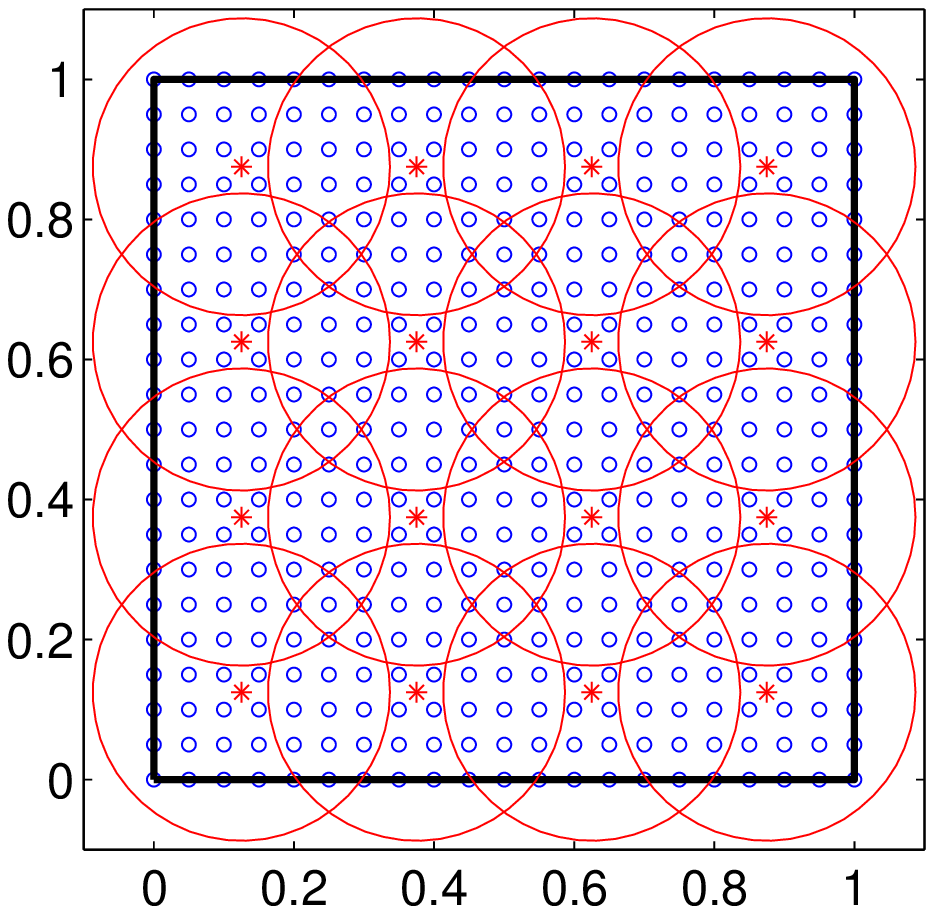}
\quad
\includegraphics[scale=0.6]{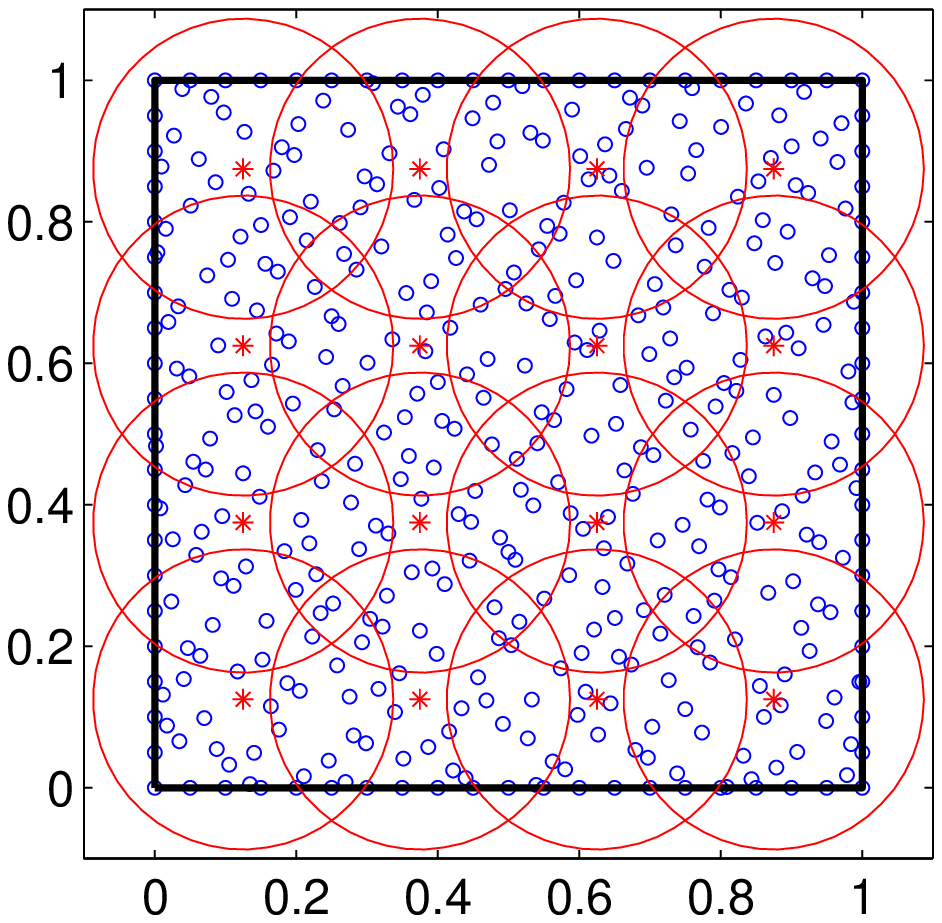}
\caption{Example of partitioning of the square domain with circular patches for uniform and Halton points (left to right).}
\label{f:4-1}
\end{figure}

In the following we focus on a wide series of experiments, which concern the 2D convection-diffusion and pseudo-parabolic problems studied in Sections \ref{cd_pp_prbs} and \ref{sec5}. In doing that, we discretize the domain $\Omega =[0, 1]\times [0, 1]$ taking some sets of two different data point distributions, which consist of a number $N$ of uniformly distributed grid nodes and quasi-random Halton points. The latter, which are a typical example of scattered node set, are generated by using the \textsc{Matlab} program \texttt{haltonseq.m} \cite{fass07}. The PU covering is composed instead of $M$ circular patches that are centred at a uniform grid of points, where the overlap of the patches is $20 \%$ of the distance between the centers \cite{saf15}. Note that, in the PUM scheme, the subdomains might be of any (regular enough) geometric shape such as circles, squares, rectangles and pentagons. However, many other possible types of patches are allowed: in fact, the only requirement for all shapes of subdomains is that they cover the domain $\Omega$. As said earlier, we here use circular subdomains so that any (possibly also mild) overlap among the different subdomains is guaranteed \cite{saf15}. An example of uniform and Halton points along with the related PU subdomains for the square domain $\Omega$ is shown in Figure \ref{f:4-1}.

We show the numerical results obtained by applying the PUM collocation method based on finite difference using some of the RBFs contained in Table \ref{tab_rbf} as local approximants. More precisely, our analysis is based on considering GA, IMQ and  M4 as basis functions in \eqref{loc_rbf}, and the compactly supported W2 as localizing function of Shepard's weight in \eqref{ShepWei}. Acting in this way, besides presenting diversified numerical tests, we compute the approximation errors, selecting \lq\lq optimal\rq\rq\ values of the shape parameter $\epsilon$. Such values of $\epsilon$ are searched for in the interval $[0.05, 12.00]$, taking as step size $0.05$. 

Notice that the accuracy of RBF-based methods highly depends upon the shape parameter $\epsilon$ of the basis functions, which is responsible for the flatness of such functions. In particular, for smooth problems the best accuracy is typically achieved when $\epsilon$ is small, but then the condition number of the linear system may become very large. In the practice, even if the PUM -- as highlighted in the next subsections -- reveals to be much better conditioned than traditional RBF techniques, the selection of shape parameters may greatly affect the accuracy of the collocation method. However, in some particularly hard situations, in which the ill-conditioning is so big to produce a very negative effect on the results, the ill-conditioning problem due to RBFs can be overtaken by using stable methods; see \cite[Chapter 12]{fass15} for an overview. Such approximation techniques allows stable computations for small values of $\epsilon$ as well.

In order to investigate the accuracy, we compute the maximum absolute error (MAE) given by
\begin{align*}
\textrm{MAE} = \max_{1\leq i \leq N}|u(\boldsymbol{x}_i,T)-p(\boldsymbol{x}_i,T)|, 
\end{align*}
where $T$ denotes the maximum time level, i.e. $T=t_{\max}=1$, whereas the stability is studied by evaluating the condition number (CN) obtained by using the \textsc{Matlab} command \texttt{condest}. Note that in this work the MAE is often referred to as \lq\lq Max Error\rq\rq. 

Although the approximation scheme described in Section \ref{cd_pp_prbs} is generally valid for any value of $\theta \in [0,1]$, in the following we focus on the case $\theta=1/2$, which identifies the famous Crank-Nicholson scheme. 

The target of our numerical analysis is thus two-fold: on the one hand, analyzing the efficiency expressed in terms of CPU times of the RBF-PUM-FD; on the other hand, verifying its accuracy and stability. In addition, we also compare our collocation method with two existing techniques: the RBF-FD method presented in \cite{esm16} and the RBF-PUM proposed in \cite{saf15}. 

Therefore, in the following sections we consider more in detail the two PDE problems: first, we show results for the convection-diffusion equation, and then for the pseudo-parabolic one.


\subsection{Results for convection-diffusion problem} \label{res1}

In this subsection we show numerical results acquired from experiments carried out for the unsteady convection-diffusion equation \eqref{e9}. Taking an appropriate initial condition and Dirichlet boundary conditions, the analytical solution is given by
\begin{equation*}
u(x,y,t)=a \exp(bt) \left( \exp(-cx)+\exp(-cy) \right),
\end{equation*}
where constants $a$ and $b$ can be chosen freely, and 
\begin{equation*}
c=\dfrac{\nu \pm \sqrt{\nu^{2}+4b\kappa}}{2\kappa}>0. 
\end{equation*}
Here, to exhibit our results after a wide experimentation with different parameters, we assume $a = 1$ and $b = 0.1$, also taking $\nu = 1$ as convection velocity and $\kappa = 1$ as diffusion coefficient.

\begin{figure}
\centering
\parbox{7.5cm}{\centering
\includegraphics[scale=0.5]{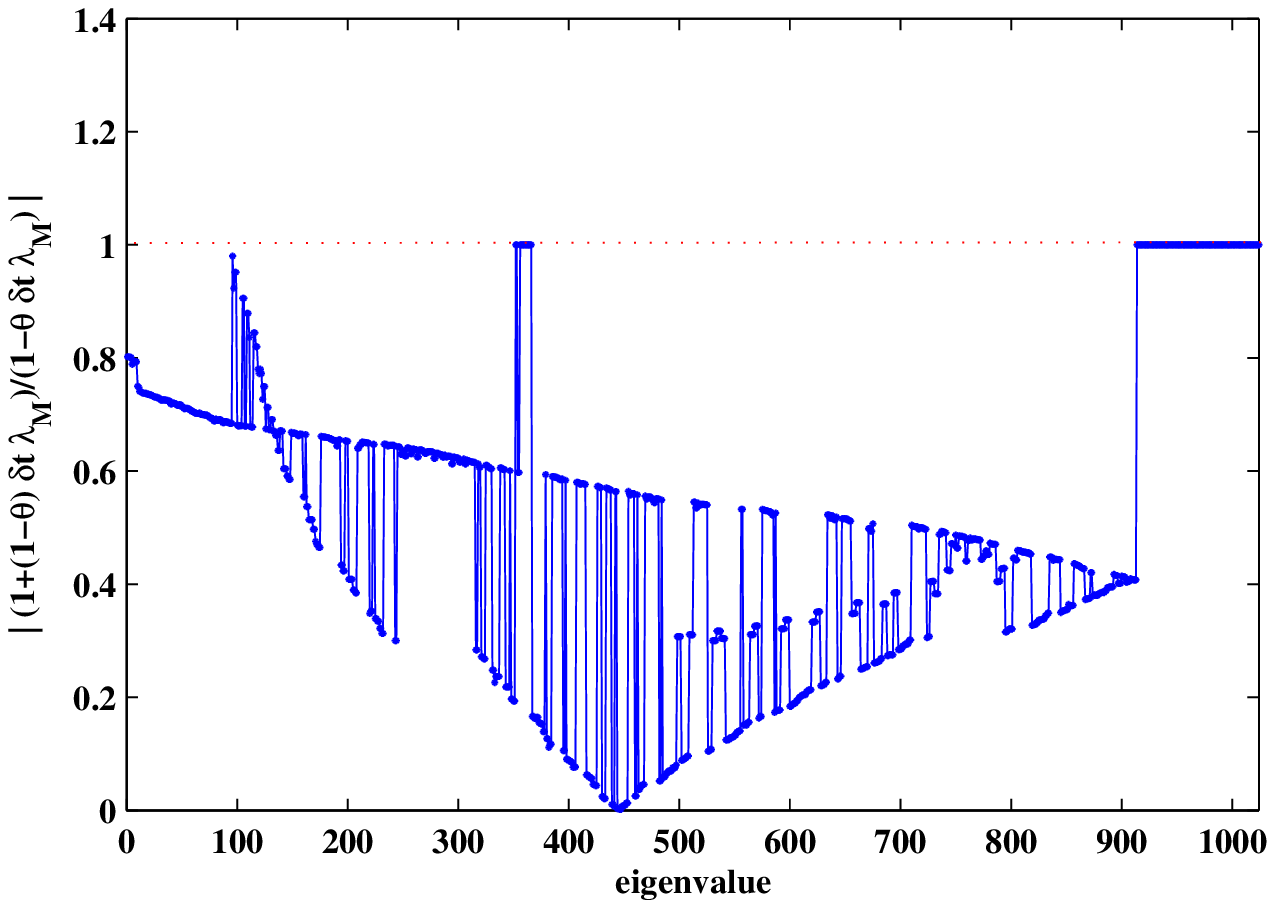}
\captionof*{figure}{Uniform points: $\epsilon = 2.85$}}
\quad
\parbox{7.5cm}{\centering
\includegraphics[scale=0.5]{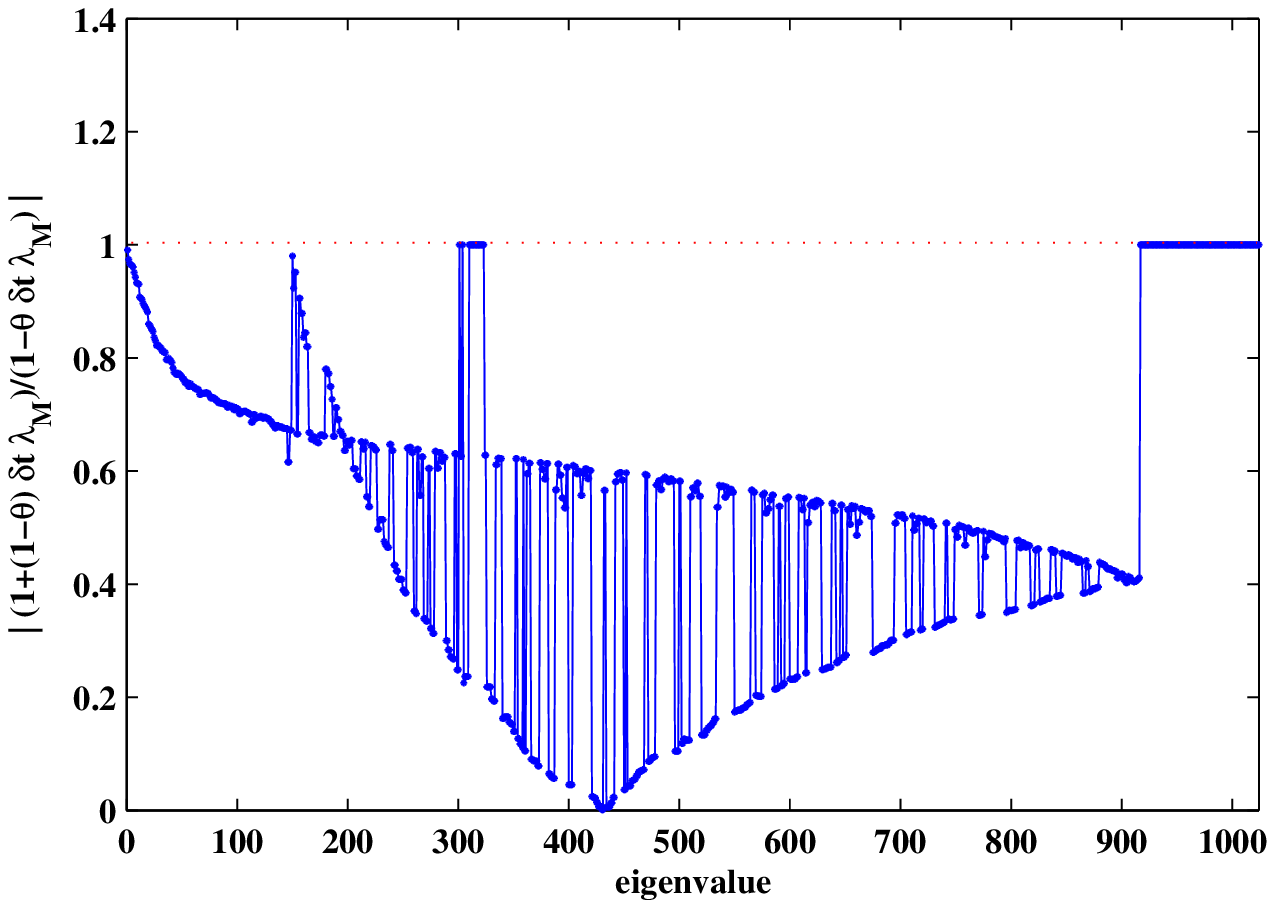}
\captionof*{figure}{Halton points: $\epsilon =2.55$}}
\caption{Numerical study of the stability of RBF-PUM-FD regarding the convection-diffusion equation with $32\times32$ uniform (left) and Halton (right) nodes and $4\times4$ patches for IMQ with $\delta t=0.001$, and $\theta=1/2$.}
\label{f1:4-1}
\end{figure}

First of all, we perform a numerical study of the stability of our method by investigating condition \eqref{e21} that depends on the eigenvalues of RBF-PUM-FD coefficient matrix for the convection-diffusion problem. As an example, in Figure \ref{f1:4-1} we highlight that condition \eqref{e21} is satisfied for all eigenvalues of matrix ${M}$ for $\theta = 1/2$. Therefore, these results also confirm from a numerical standpoint that our collocation method is unconditionally stable.

\begin{figure}
\centering
\parbox{5cm}{\centering
\includegraphics[scale=0.4]{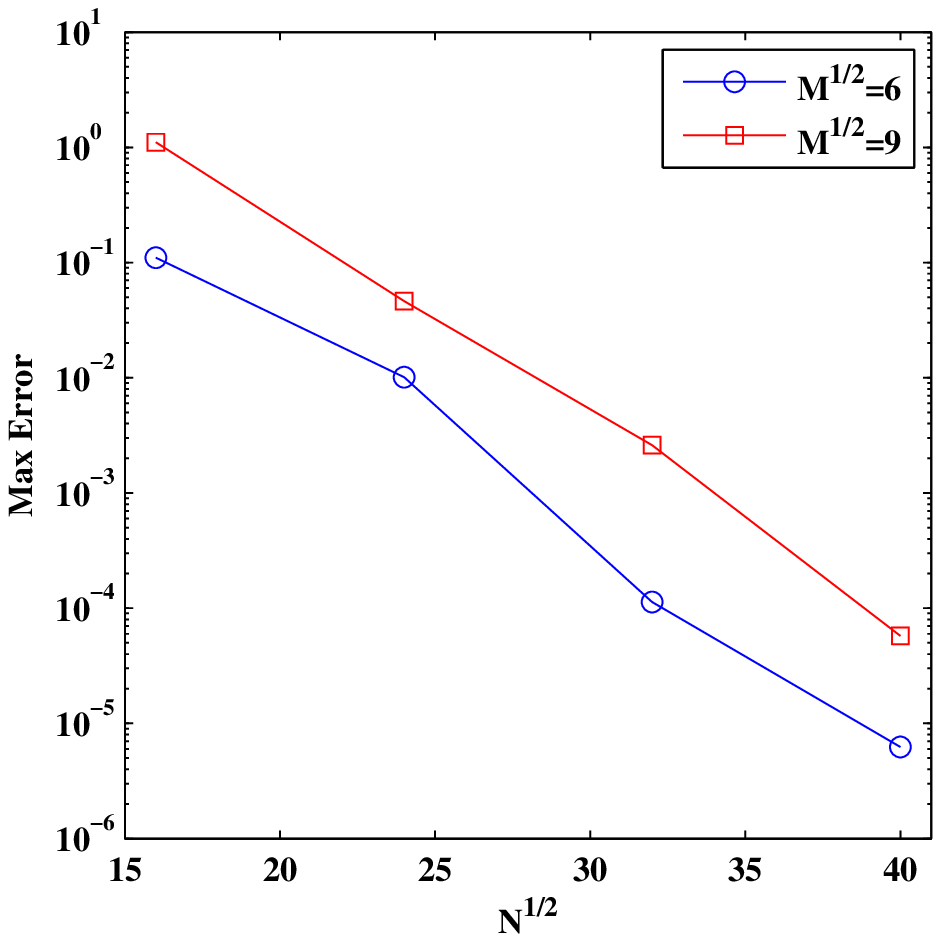}
\captionof*{figure}{Uniform points: IMQ, $\epsilon = 3.20$}}
\quad
\parbox{5cm}{\centering
\includegraphics[scale=0.4]{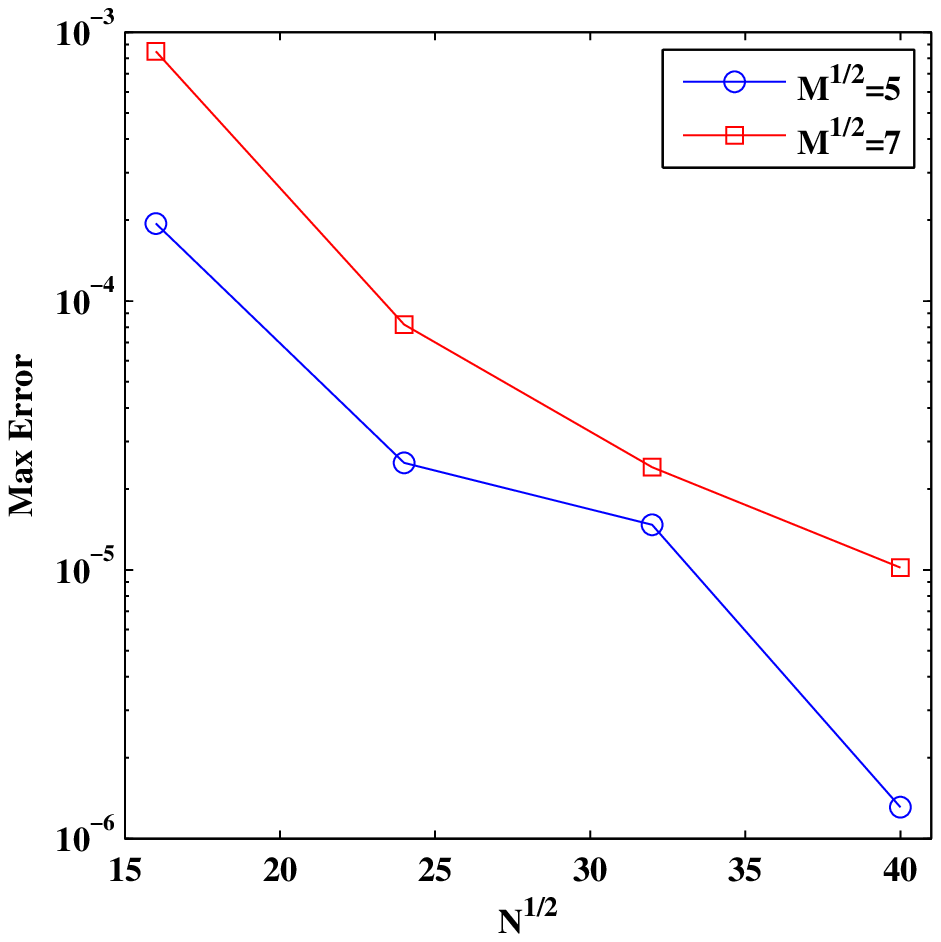}
\captionof*{figure}{Uniform points: M4, $\epsilon = 0.20$}}
\quad
\parbox{5cm}{\centering
\includegraphics[scale=0.4]{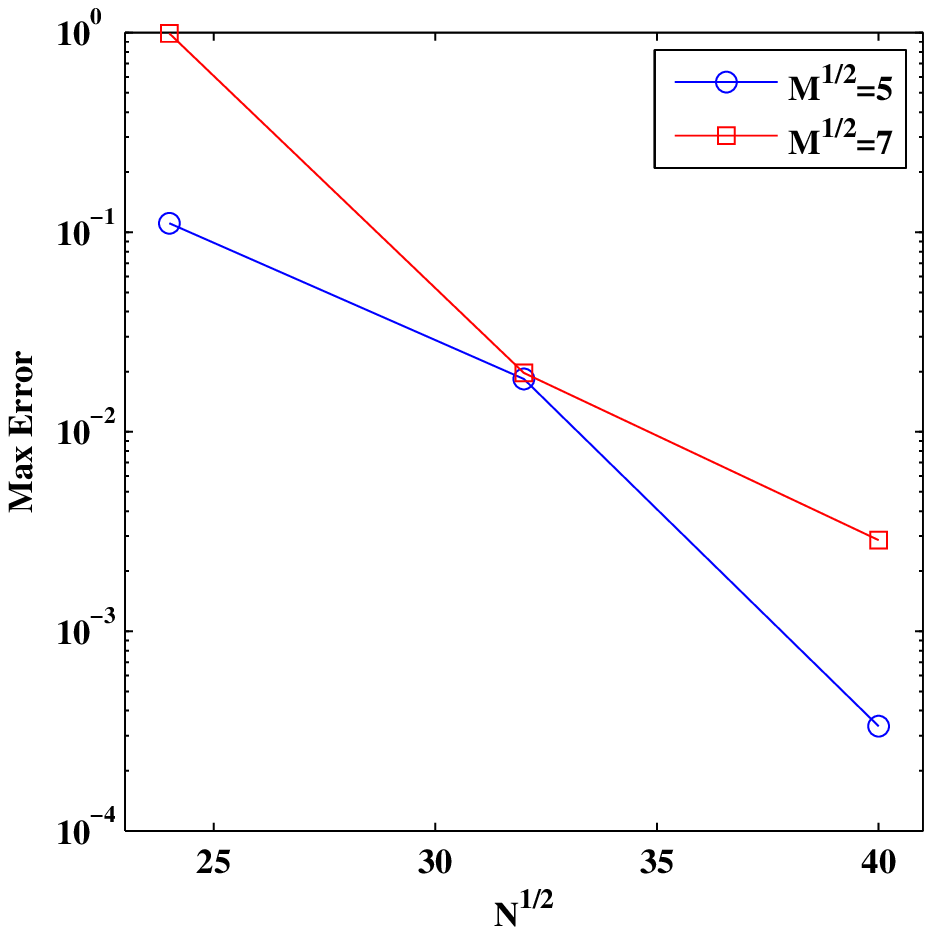}
\captionof*{figure}{Uniform points: GA, $\epsilon = 9.20$}}\\
\parbox{5cm}{\centering
\includegraphics[scale=0.4]{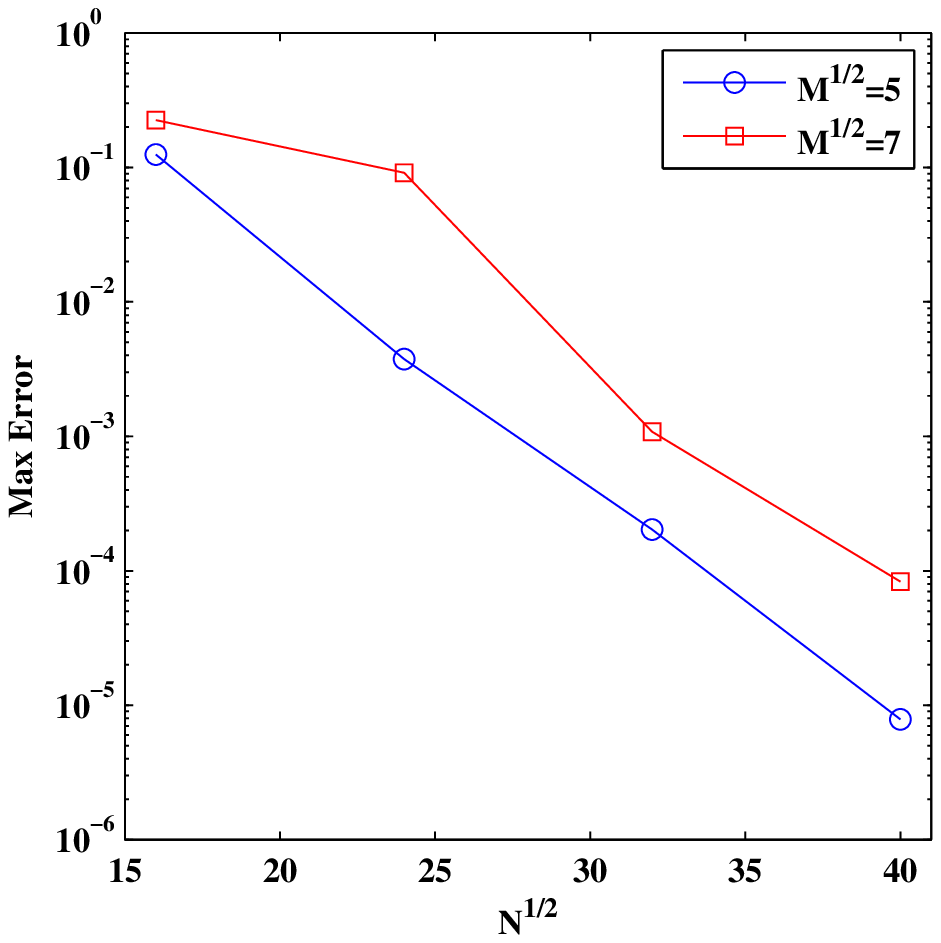}
\captionof*{figure}{Halton points: IMQ, $\epsilon = 4.25$}}
\quad
\parbox{5cm}{\centering
\includegraphics[scale=0.4]{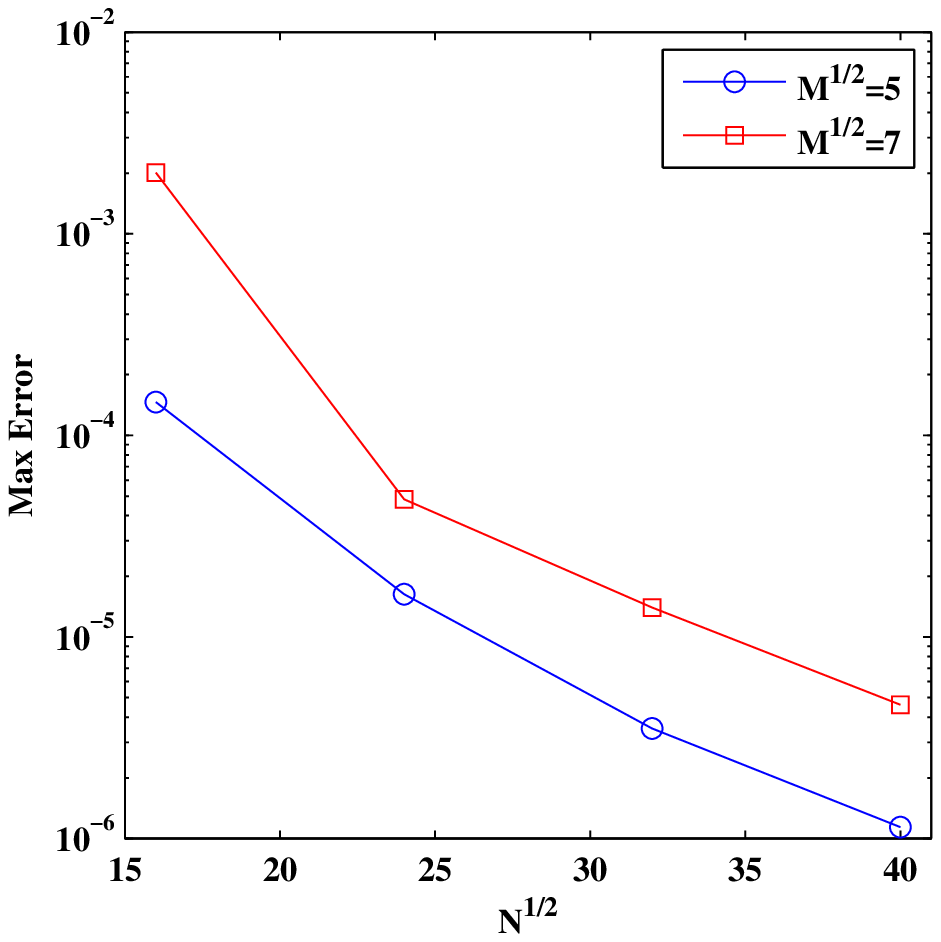}
\captionof*{figure}{Halton points: M4, $\epsilon = 0.20$}}
\quad
\parbox{5cm}{\centering
\includegraphics[scale=0.4]{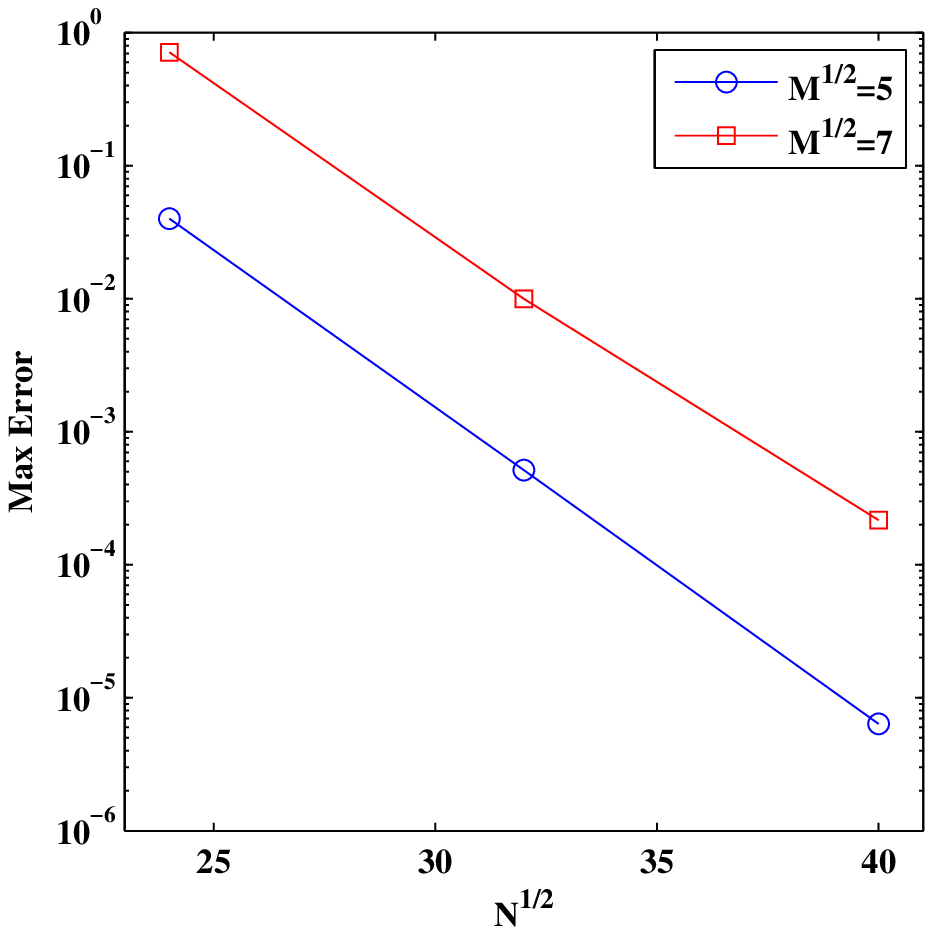}
\captionof*{figure}{Halton points: GA, $\epsilon = 7.65$}}
\caption{Error in the convection-diffusion equation against the problem size with respect to the number of partitions using $\delta t=0.001$ and $T=1$.}\label{f:2}
\end{figure}

In Figure \ref{f:2} we analyze the error behavior of our RBF-PUM-FD collocation method by varying the number of points. In particular, we show how the error changes with respect to the number of partitions (or subdomains) using IMQ, GA and M4 in the RBF-PUM-FD scheme for uniform and Halton points. From this study, we can note that the smaller (larger) the number of partitions is the smaller (larger) the error becomes. Accordingly, increasing (decreasing) the number of points per partition results in a better accuracy, but the linear system \eqref{e15} becomes denser (sparser) thus requiring more (less) time for its solution. In our experiments we carefully analyzed this issue and in the following we report results in which we achieved good trade-off between accuracy and computational efficiency.  

\begin{figure}
\centering
\includegraphics[scale=0.6]{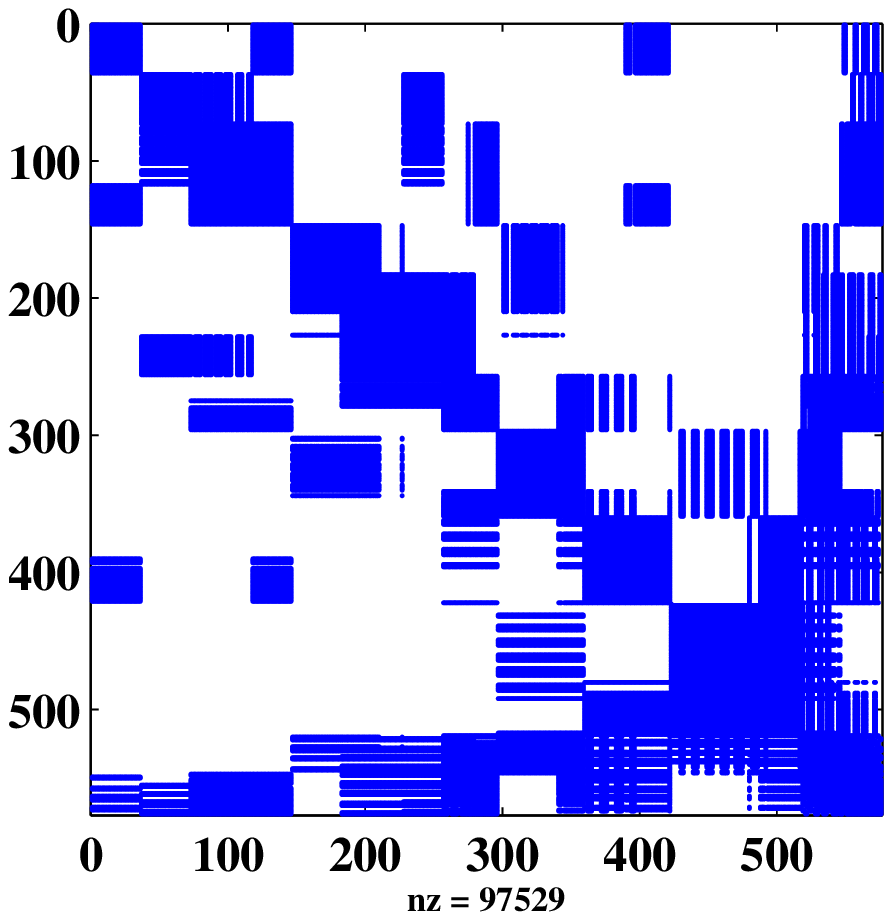}
\hskip -1.5cm
\includegraphics[scale=0.6]{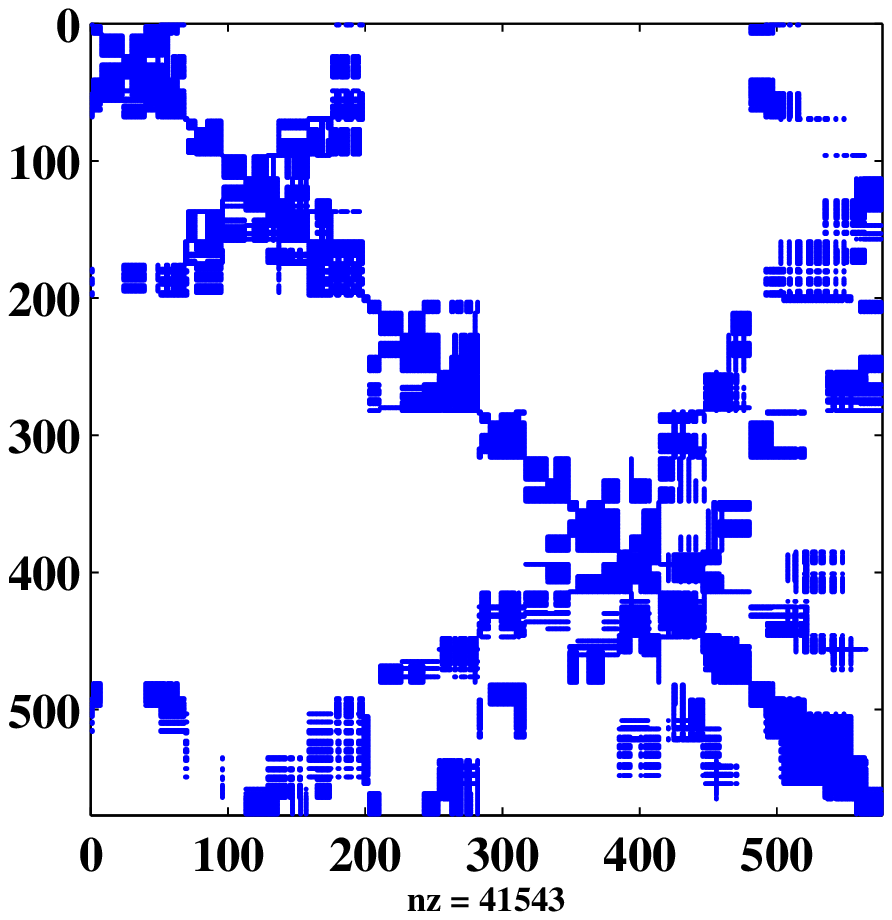}
\caption{Sparsity structure of coefficient matrix regarding the convection-diffusion equation with $24\times24$ uniform nodes and $3\times3$ patches (left) and $5\times5$ patches (right) for IMQ with $\epsilon=1.80$.}\label{f1:1}
\end{figure}

One of the main advantages of the RBF-PUM-FD is to enjoy a great flexibility. In fact, this scheme allows to choose quite freely the number of partitions to cover the domain. Generally, a covering with smaller partitions leads to worse approximation results but, at the same time, it turns out to be less expensive from the computational standpoint. This is essentially due to the greater sparsity of the linear system. In order to stress sparsity of the coefficient matrix, in Figure \ref{f1:1} we show the sparse structure generated by the RBF-PUM-FD collocation scheme for the convection-diffusion equation considering IMQ and two different domain partitions. As it is evident from these graphs, fixed the number of nodes more patches lead to more sparsity: this increase contributes to a better efficiency of the numerical scheme but at the same time lessens its precision. Even if only the diffusion term is present, the matrices are non-symmetric due to the collocation involving the partition of unity weight functions.

\begin{figure}
\centering
\parbox{7.5cm}{\centering
\includegraphics[scale=0.325]{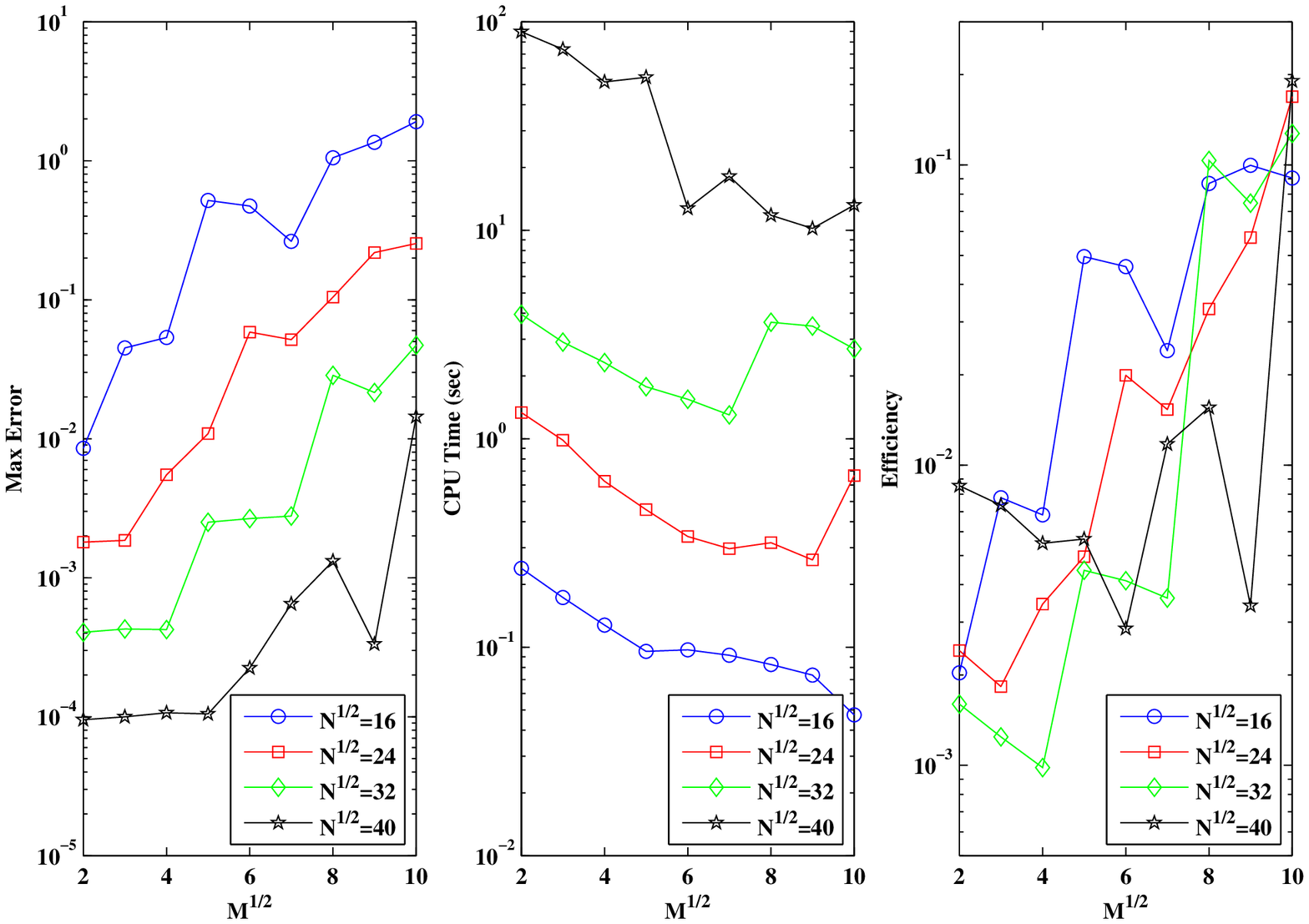}
\captionof*{figure}{Uniform points: IMQ, $\epsilon = 5$}}
\parbox{7.5cm}{\centering
\includegraphics[scale=0.325]{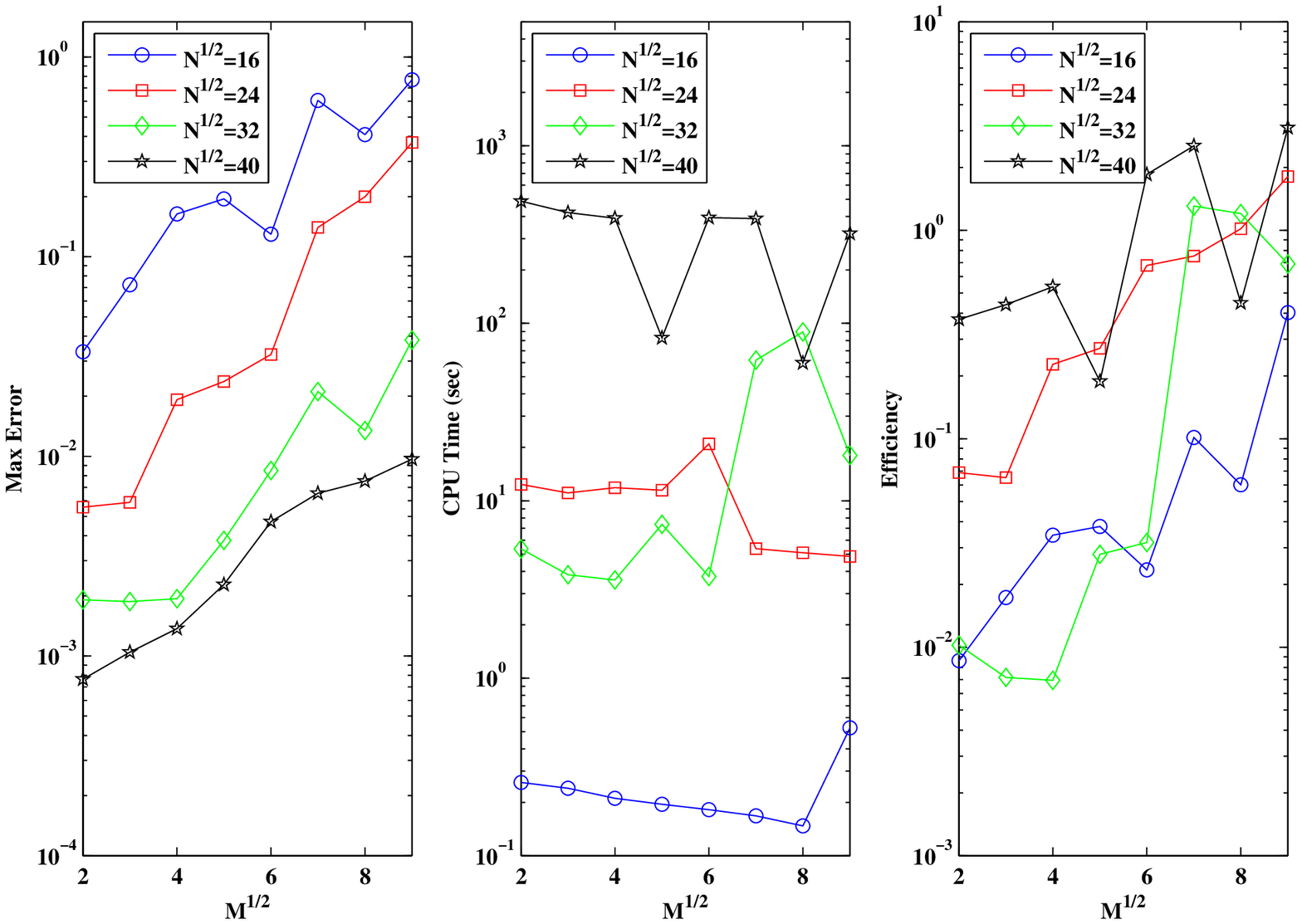}
\captionof*{figure}{Halton points: M4, $\epsilon = 10$}}
\caption{Left: Error in the convection-diffusion equation against the number of partitions in one spatial dimension. Centre: Computational time against the number of partitions in one spatial dimension. Right: Efficiency computed as product between the error and CPU time. Parameters: $\delta t=0.001$ and $T=1$.}\label{f:1}
\end{figure}

Additionally, in Figure \ref{f:1} we show the error versus the number of partitions in one spatial dimension (left), the associated CPU time (centre), and the computational efficiency as a product of the two (right). Although these tests have been carried out for various RBFs providing similar results, for shortness we report only a single case for uniform points using IMQ and one for Halton points using M4. These graphs give us some useful information on the possible choice of parameters; for instance, efficiency can provide us a suggestion to select a good number of partitions compared to MAE and CPU time.

\begin{figure}
\centering
\includegraphics[scale=0.5]{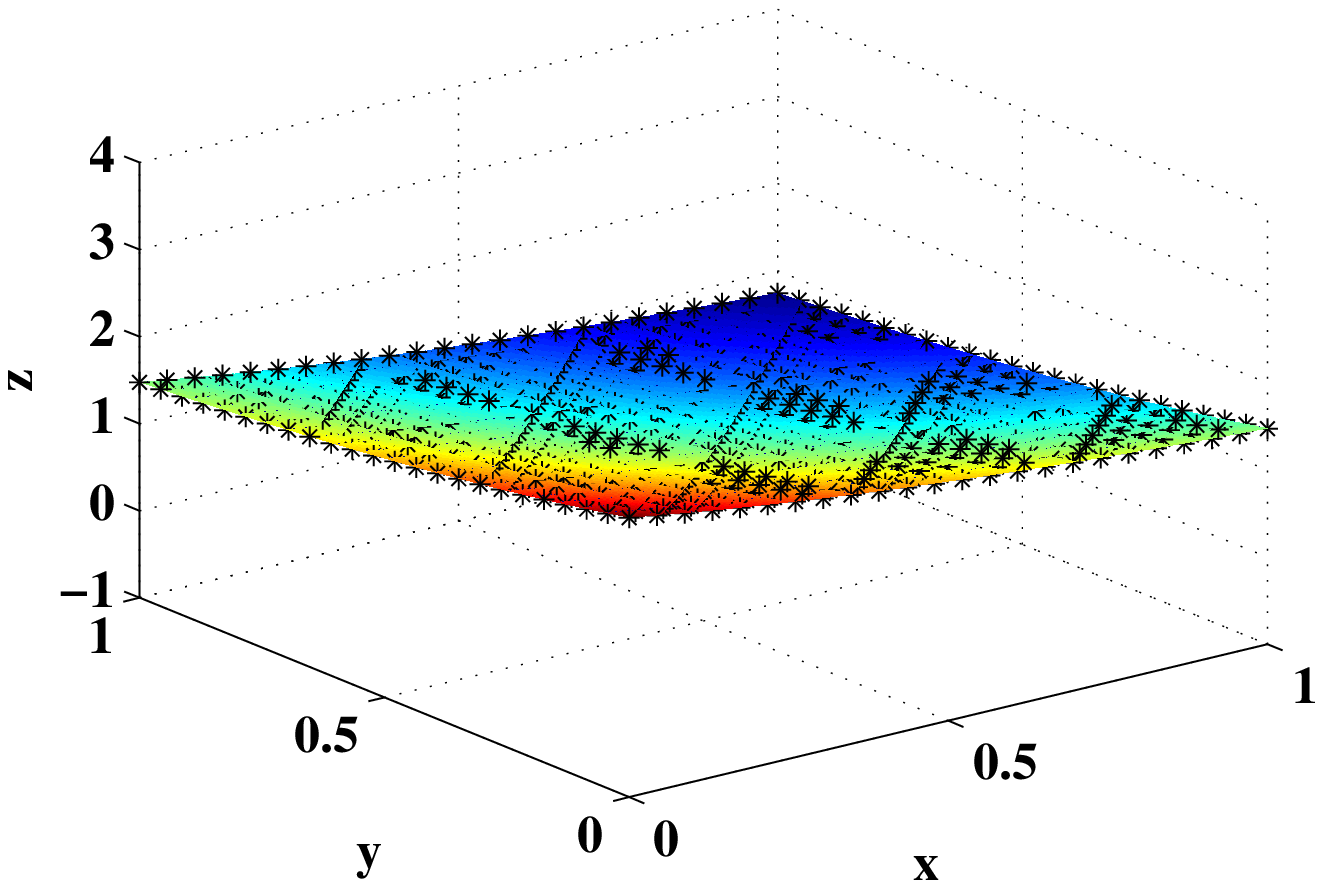}
\quad
\includegraphics[scale=0.5]{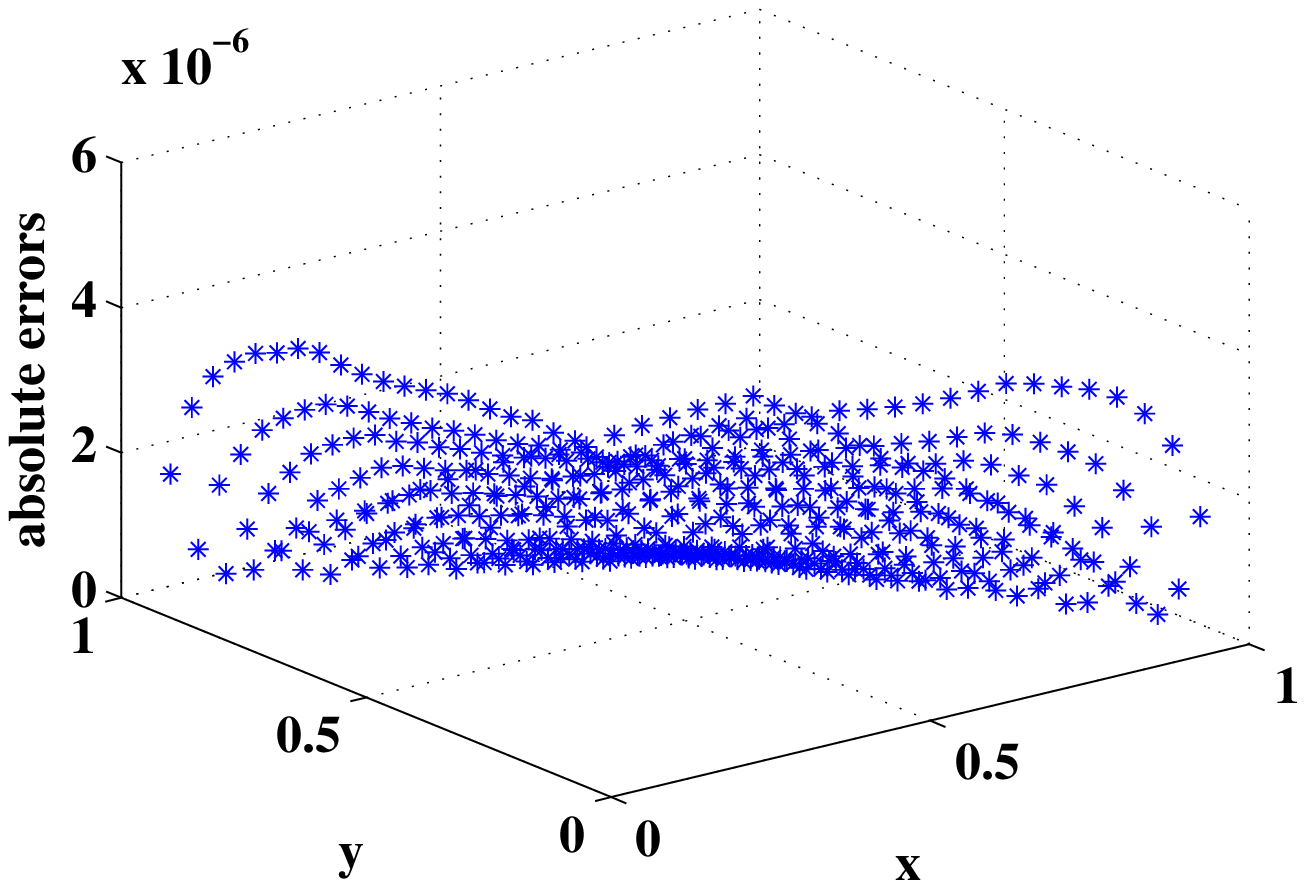}
\caption{Surface of exact solution \eqref{e9} and its approximation given per points (left) along with the related absolute errors (right) on $24\times24$ uniform nodes and $3\times3$ patches, using IMQ for $\epsilon = 1.80$ with $\delta t=0.001$ and $T=1$.}
\label{f1::1}
\end{figure}

In Tables \ref{tab:1}-\ref{tab:t} we present an extensive and detailed analysis on the RBF-PUM-FD collocation method, which is also compared with the global RBF-FD method studied in \cite{esm16}. More precisely, considering four sets of uniform and Halton points and IMQ as local RBF approximant, we report the MAE computed with \lq\lq optimal\rq\rq\ $\epsilon$, the CPU time required to solve the sparse linear system \eqref{e15} and the CN. Furthermore, as an example, using a set of uniform nodes and IMQ as local approximant in the PUM scheme, in Figure \ref{f1::1} we graphically represent surface of exact solution and its approximation given per points (left); on the right, instead, we show the related absolute errors computed at the same collocation points. In all cases we fix $\delta t = 0.001$ and $T=1$. From this study we can note a quite uniform behavior: on the one hand, accuracy of the two methods (RBF-PUM-FD vs RBF-FD) is almost similar, with usually a slight advantage for the RBF-PUM-FD method, but on the other hand we observe a significant reduction of CPU time and CN for the RBF-PUM-FD method compared to the RBF-FD one. Note that, in general, for the RBF-PUM-FD scheme we found a good compromise between accuracy and efficiency assuming that the ratio $\sqrt{N}/\sqrt{M}=8$. Obviously, the flexibility of our method enables also to consider different values as it can occur in some situations.

\begin{table}
\begin{center} 
{\small
\begin{tabular}{ccccccccccc} \toprule
 & & & & & RBF-PUM-FD & $\quad$ & & & & RBF-FD \\ \cmidrule{2-6} \cmidrule{8-11}
$\sqrt{N}$ & $\sqrt{M}$& $\epsilon$ & MAE & time & CN & $\quad$ &$\epsilon$ & MAE & time & CN \\ \midrule
$16$ & $2$& $1.35$ & $1.91 \times 10^{-5}$ & $0.22$ & $4.08\times 10^{+01}$ & $\quad$ & $1.55$ & $3.85 \times 10^{-5}$ & $6.46$ & $3.89 \times 10^{+17}$\\ \midrule
$24$ & $3$& $1.80$ & $4.51 \times 10^{-6}$ & $5.88$ & $1.28\times 10^{+04}$ & $\quad$ & $2.40$ & $2.89 \times 10^{-5}$ & $21.06$ & $5.69 \times 10^{+17}$\\ \midrule
$32$ & $4$& $2.85$ & $1.19 \times 10^{-5}$ & $12.50$ & $2.80\times 10^{+02}$ & $\quad$ & $3.25$ & $2.49 \times 10^{-5}$ & $130.88$ & $4.43 \times 10^{+17}$\\ \midrule
$40$ & $5$& $3.20$ & $6.22 \times 10^{-6}$ & $21.45$ & $3.54\times 10^{+02}$ & $\quad$ & $4.25$ & $2.95 \times 10^{-5}$ & $490.51$ & $2.91 \times 10^{+17}$\\ \bottomrule
\end{tabular}
}
\caption{\small Comparison between RBF-PUM-FD and RBF-FD \cite{esm16} using IMQ for the convection-diffusion equation: MAE, CPU time (in seconds) and CN for uniform points with $\delta t=0.001$ and $T=1$.} 
\label{tab:1}
\end{center}
\end{table}

\begin{table}
\begin{center} 
{\small
\begin{tabular}{ccccccccccc} \toprule
 & & & & & RBF-PUM-FD & $\quad$ & & & & RBF-FD \\ \cmidrule{2-6} \cmidrule{8-11}
$\sqrt{N}$ & $\sqrt{M}$& $\epsilon$ & MAE & time & CN & $\quad$ &$\epsilon$ & MAE & time & CN \\ \midrule
$16$ & $2$& $1.10$ & $1.16 \times 10^{-6}$ & $1.87$ & $8.48\times 10^{+02}$ & $\quad$ & $1.30$ & $1.25 \times 10^{-6}$ & $6.54$ & $7.24 \times 10^{+17}$\\ \midrule
$24$ & $3$& $1.85$ & $9.51 \times 10^{-7}$ & $15.95$ & $1.19\times 10^{+05}$ & $\quad$ & $2.25$ & $6.87 \times 10^{-7}$ & $20.19$ & $6.07 \times 10^{+17}$\\ \midrule
$32$ & $4$& $2.55$ & $5.22 \times 10^{-7}$ & $32.27$ & $3.95\times 10^{+04}$ & $\quad$ & $3.05$ & $1.51 \times 10^{-6}$ & $151.95$ & $1.09 \times 10^{+18}$\\ \midrule
$40$ & $5$& $3.25$ & $5.92 \times 10^{-7}$ & $84.51$ & $3.02\times 10^{+04}$ & $\quad$ & $3.85$ & $5.11 \times 10^{-7}$ & $485.76$ & $1.03 \times 10^{+18}$\\ \bottomrule
\end{tabular}
}
\caption{\small Comparison between RBF-PUM-FD and RBF-FD \cite{esm16} using IMQ for the convection-diffusion equation: MAE, CPU time (in seconds) and CN for Halton points with $\delta t=0.001$ and $T=1$.}  
\label{tab:1-h}
\end{center}
\end{table}

Moreover, since in previous tables execution time refers to the one needed to solve the sparse system \eqref{e15}, in Table \ref{tab:t} we focus on total CPU time for both RBF-PUM-FD and RBF-FD methods. From these results, we can notice as the RBF-PUM-FD needs a setup cost due to data organization in the PUM scheme, while the time in the global RBF-FD is roughly the same (cf. Tables \ref{tab:1}-\ref{tab:1-h} with Table \ref{tab:t}).

\begin{table}
\begin{center} 
{\small
\begin{tabular}{cccccccc} \toprule
 & & Uniform points & & $\quad$ & & Halton points & \\ \cmidrule{2-4} \cmidrule{6-8}
$\sqrt{N}$ & $\sqrt{M}$ & RBF-PUM-FD & RBF-FD & $\quad$ & $\sqrt{M}$ & RBF-PUM-FD & RBF-FD \\ \midrule
$16$ & $2$ & $1.36$ & $6.89$ & $\quad$ & $2$ & $3.04$ & $6.94$ \\ \midrule
$24$ & $3$ & $10.56$ & $21.52$ & $\quad$ & $3$ & $19.33$ & $21.94$ \\ \midrule
$32$ & $4$ & $27.17$ & $131.19$ & $\quad$ & $4$ & $47.61$ & $152.11$ \\ \midrule
$40$ & $5$ & $44.69$ & $491.60$ & $\quad$ & $5$ & $125.08$ & $489.94$ \\ \bottomrule
\end{tabular}
}
\caption{\small Comparison of total CPU time (in seconds) between RBF-PUM-FD and RBF-FD \cite{esm16} using IMQ for the convection-diffusion equation with $\delta t=0.001$ and $T=1$.} 
\label{tab:t}
\end{center}
\end{table}

Finally, in order to give a more complete analysis of our collocation method (that relies on combination of finite difference discretization in time domain and local meshfree RBF-based discretization in space), in Tables \ref{tab:t1-m}-\ref{tab:t1-m-h} we compare accuracy (MAE), efficiency (time) and stability (CN) of the RBF-PUM-FD method with the RBF-PUM scheme proposed in \cite{saf15} using M4. This study points out essentially that our method is generally accurate as the RBF-PUM, more expensive than the RBF-PUM from a computational viewpoint, but also less ill-conditioned.

\begin{table}
\begin{center} 
{\small
\begin{tabular}{ccccccccccc} \toprule
 & & & & & RBF-PUM-FD & $\quad$ & & & & RBF-PUM  \\ \cmidrule{3-6} \cmidrule{8-11}
$\sqrt{N}$ & $\sqrt{M}$& $\epsilon$ & MAE & time & CN & $\quad$ &$\epsilon$ & MAE & time & CN \\ \midrule
$16$ & $2$& $0.05$ & $1.83 \times 10^{-5}$ & $1.44$ & $3.81\times 10^{+01}$ & $\quad$ & $0.05$ & $1.83 \times 10^{-5}$ & $1.85$ & $1.42 \times 10^{+03}$\\ \midrule
$24$ & $3$& $0.10$ & $5.37 \times 10^{-6}$ & $6.19$ & $1.46\times 10^{+01}$ & $\quad$ & $0.10$ & $5.37 \times 10^{-6}$ & $6.39$ & $1.23 \times 10^{+03}$\\ \midrule
$32$ & $4$& $0.15$ & $2.36 \times 10^{-6}$ & $21.96$ & $2.40\times 10^{+01}$ & $\quad$ & $0.10$ & $2.19 \times 10^{-6}$ & $36.60$ & $2.48 \times 10^{+04}$\\ \midrule
$40$ & $5$& $0.20$ & $1.31 \times 10^{-6}$ & $49.13$ & $3.80\times 10^{+01}$ & $\quad$ & $0.15$ & $1.19 \times 10^{-6}$ & $58.21$ & $7.94 \times 10^{+03}$\\ \bottomrule
\end{tabular}
}
\caption{\small Comparison between RBF-PUM-FD and RBF-PUM \cite{saf15} using M4 for the convection-diffusion equation: MAE, CPU time (in seconds) and CN for uniform points with $\delta t=0.001$ and $T=1$.} 
\label{tab:t1-m}
\end{center}
\end{table}

\begin{table}
\begin{center} 
{\small
\begin{tabular}{ccccccccccc} \toprule
 & & & & & RBF-PUM-FD & $\quad$ & & & & RBF-PUM  \\ \cmidrule{3-6} \cmidrule{8-11}
$\sqrt{N}$ & $\sqrt{M}$& $\epsilon$ & MAE & time & CN & $\quad$ &$\epsilon$ & MAE & time & CN \\ \midrule
$16$ & $2$& $0.10$ & $2.82 \times 10^{-5}$ & $1.46$ & $4.95\times 10^{+01}$ & $\quad$ & $0.10$ & $2.82 \times 10^{-5}$ & $1.85$ & $4.10 \times 10^{+03}$\\ \midrule
$24$ & $3$& $0.10$ & $5.81 \times 10^{-6}$ & $18.32$ & $8.81\times 10^{+02}$ & $\quad$ & $0.10$ & $5.82 \times 10^{-6}$ & $6.14$ & $3.02 \times 10^{+04}$\\ \midrule
$32$ & $4$& $0.15$ & $1.88 \times 10^{-6}$ & $38.47$ & $2.01\times 10^{+03}$ & $\quad$ & $0.15$ & $1.86 \times 10^{-6}$ & $18.24$ & $1.29 \times 10^{+05}$\\ \midrule
$40$ & $5$& $0.20$ & $1.14 \times 10^{-6}$ & $121.19$ & $1.51\times 10^{+03}$ & $\quad$ & $0.15$ & $1.09 \times 10^{-6}$ & $63.71$ & $5.97 \times 10^{+05}$\\ \bottomrule
\end{tabular}
}
\caption{\small Comparison between RBF-PUM-FD and RBF-PUM \cite{saf15} using M4 for the convection-diffusion equation: MAE, CPU time (in seconds) and CN for Halton points with $\delta t=0.001$ and $T=1$.} 
\label{tab:t1-m-h}
\end{center}
\end{table}


\subsection{Results for pseudo-parabolic problem} \label{res2}
In this subsection we present some of the numerical experiments we made to test our method for the pseudo-\-pa\-ra\-bo\-lic equation \eqref{e9-1}. Taken a finite domain $\Omega=\lbrace (x,y,t):0<x<1,0<y<1,0\leq t \leq T \rbrace$, we consider Eqs. \eqref{e9-1}-\eqref{e14-1} with initial condition
\begin{equation*}
u(x,y,0)=\cos(x)+\sin(y),\quad (x,y)\in [0,1]^{2}.
\end{equation*}
The exact solution to the pseudo-parabolic equation is
\begin{equation*}
u(x,y,t)=\exp(2t)(\cos(x)+\sin(y)),
\end{equation*}
while the function $f$ that appears in \eqref{e9-1} is given by
\begin{equation*}
f(x,y,t)=(2+\alpha+2\beta)\exp(2t)(\cos(x)+\sin(y)).
\end{equation*}
The experiments we report have been performed assuming as parameters the values of $\alpha = 1$ and $\beta = 0.00025$.

\begin{figure}
\centering
\parbox{7.5cm}{\centering
\includegraphics[scale=0.5]{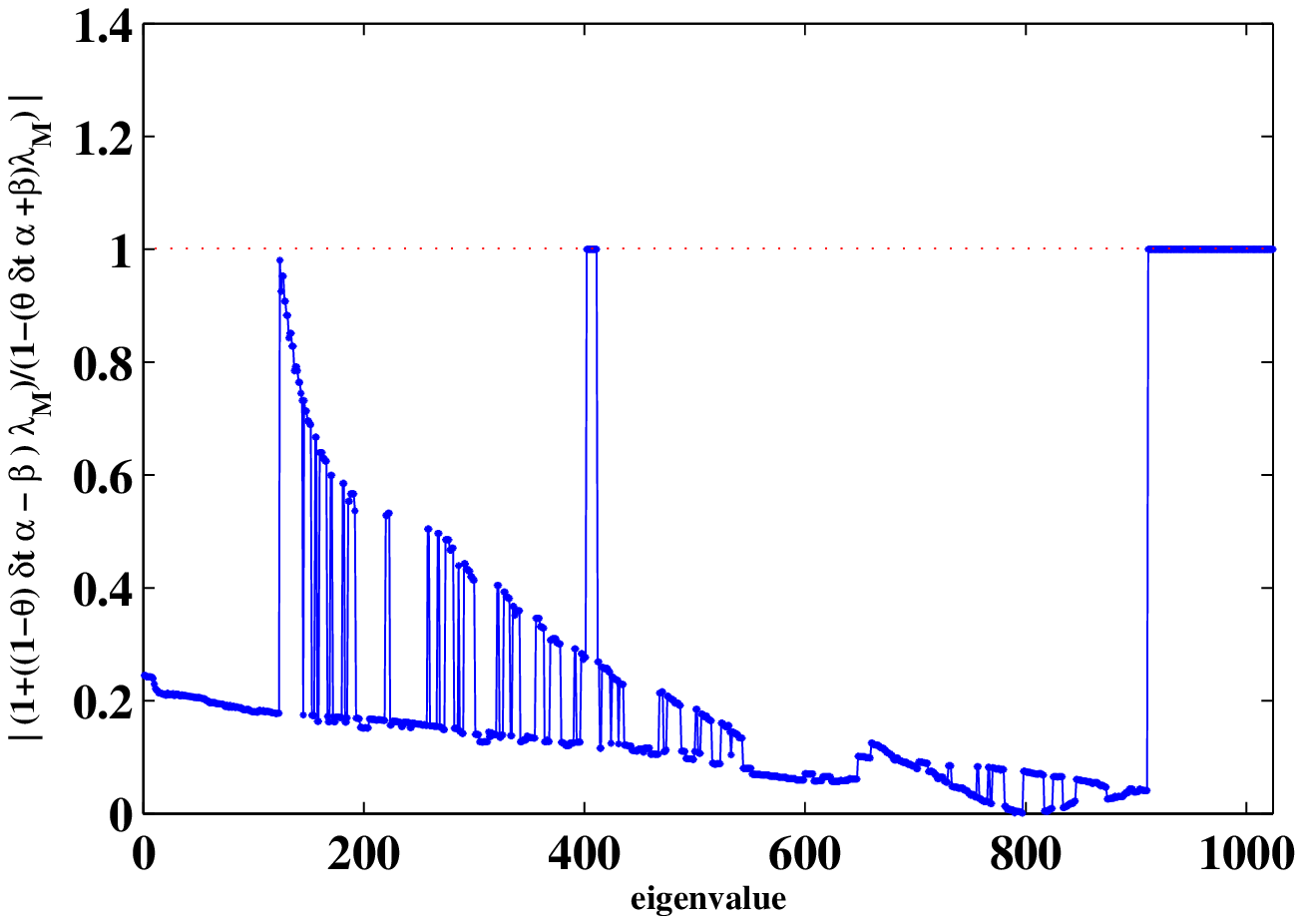}
\captionof*{figure}{Uniform points: $\epsilon = 7.60$}}
\quad
\parbox{7.5cm}{\centering
\includegraphics[scale=0.5]{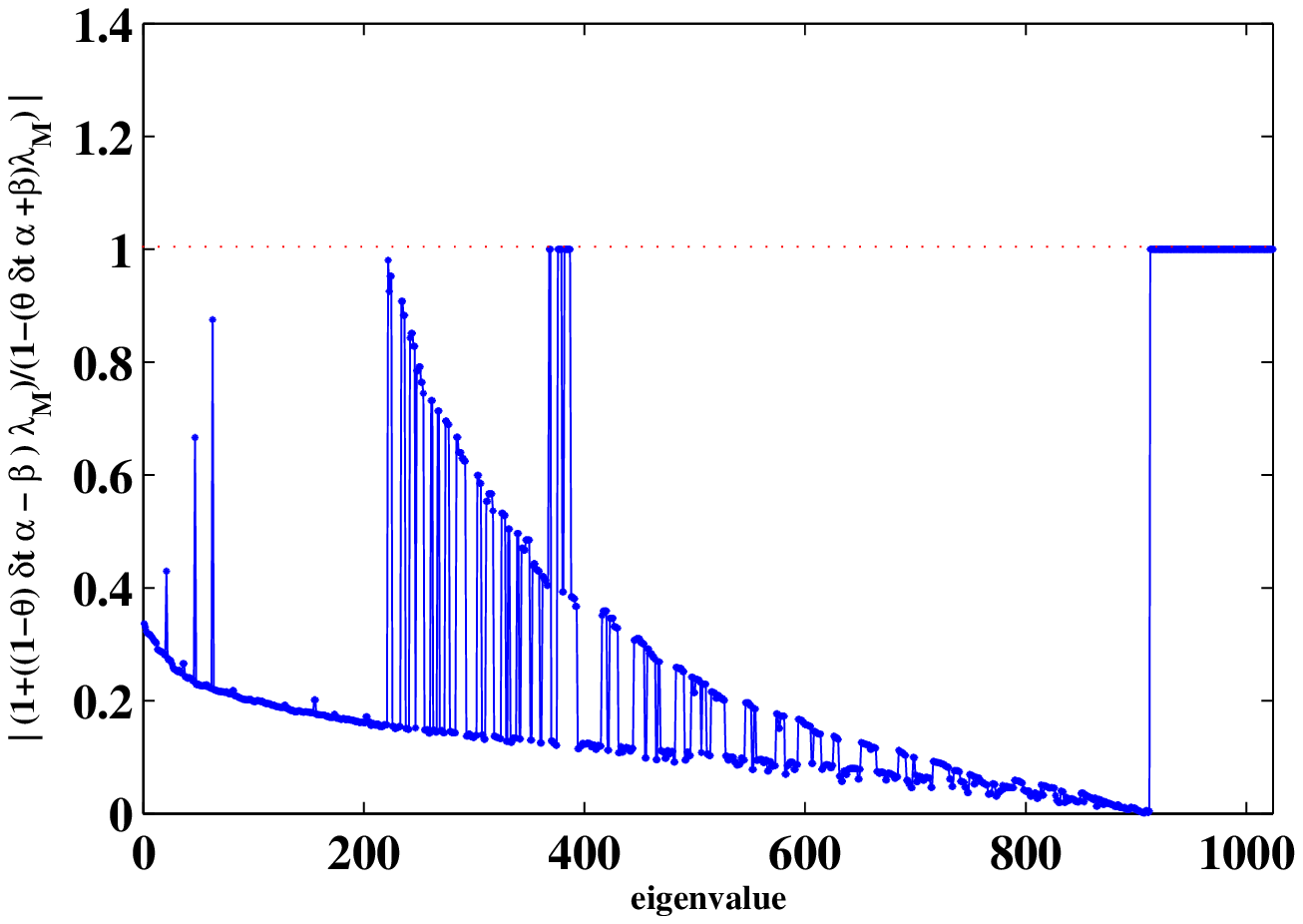}
\captionof*{figure}{Halton points: $\epsilon =5.50$}}
\caption{Numerical study of the stability of RBF-PUM-FD regarding the pseudo-parabolic equation with $32\times32$ uniform (left) and Halton (right) nodes and $4\times4$ patches for GA with $\delta t=0.001$, and $\theta=1/2$.}
\label{f1b:4-1}
\end{figure}

At first, we analyze numerically the stability of our RBF-PUM-FD scheme by verifying validity of the inequality \eqref{e21-1} for the pseudo-parabolic equation. In particular, in Figure \ref{f1b:4-1} we show that condition \eqref{e21-1} is satisfied for all eigenvalues of matrix ${M}$ for $\theta = 1/2$. These graphs point out that our collocation scheme turns out to be unconditionally stable.

\begin{figure}
\centering
\parbox{5cm}{\centering
\includegraphics[scale=0.37]{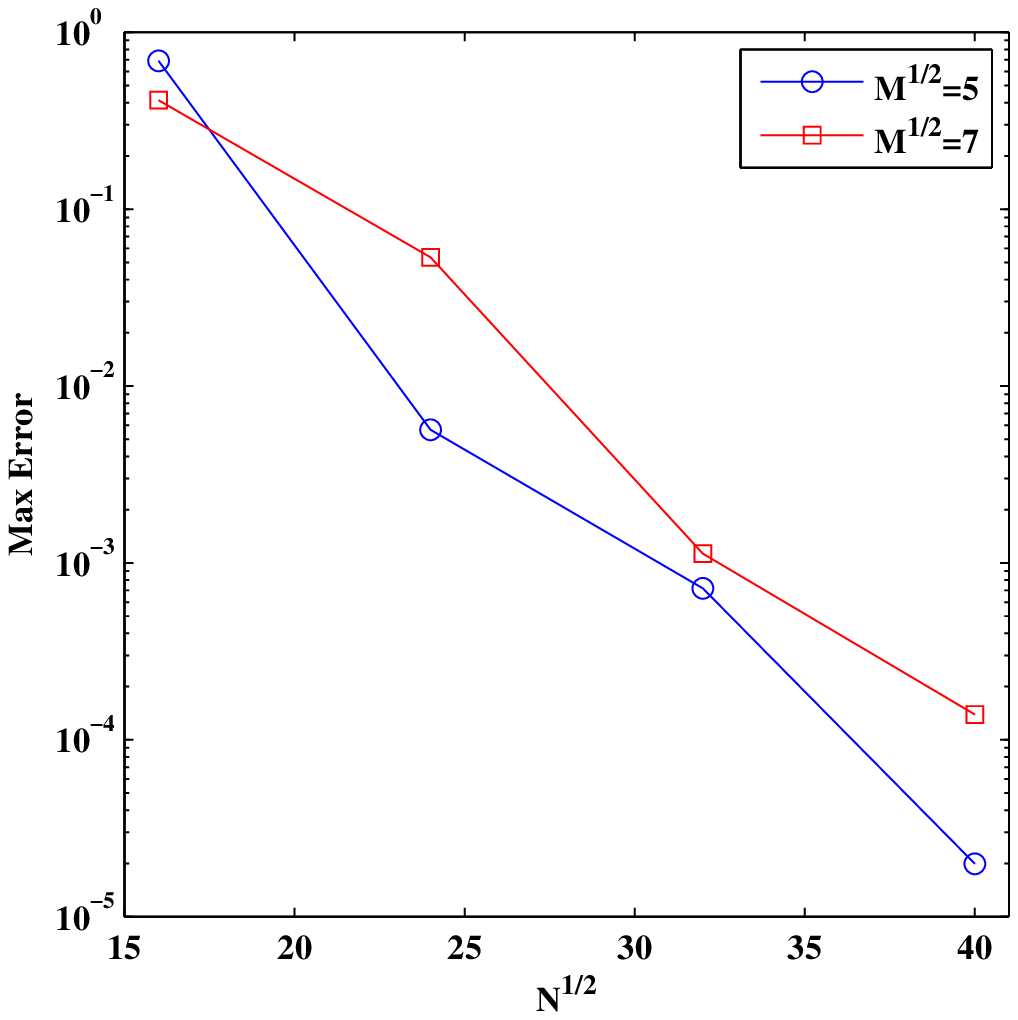}
\captionof*{figure}{Uniform points: IMQ, $\epsilon = 3.20$}}
\quad
\parbox{5cm}{\centering
\includegraphics[scale=0.37]{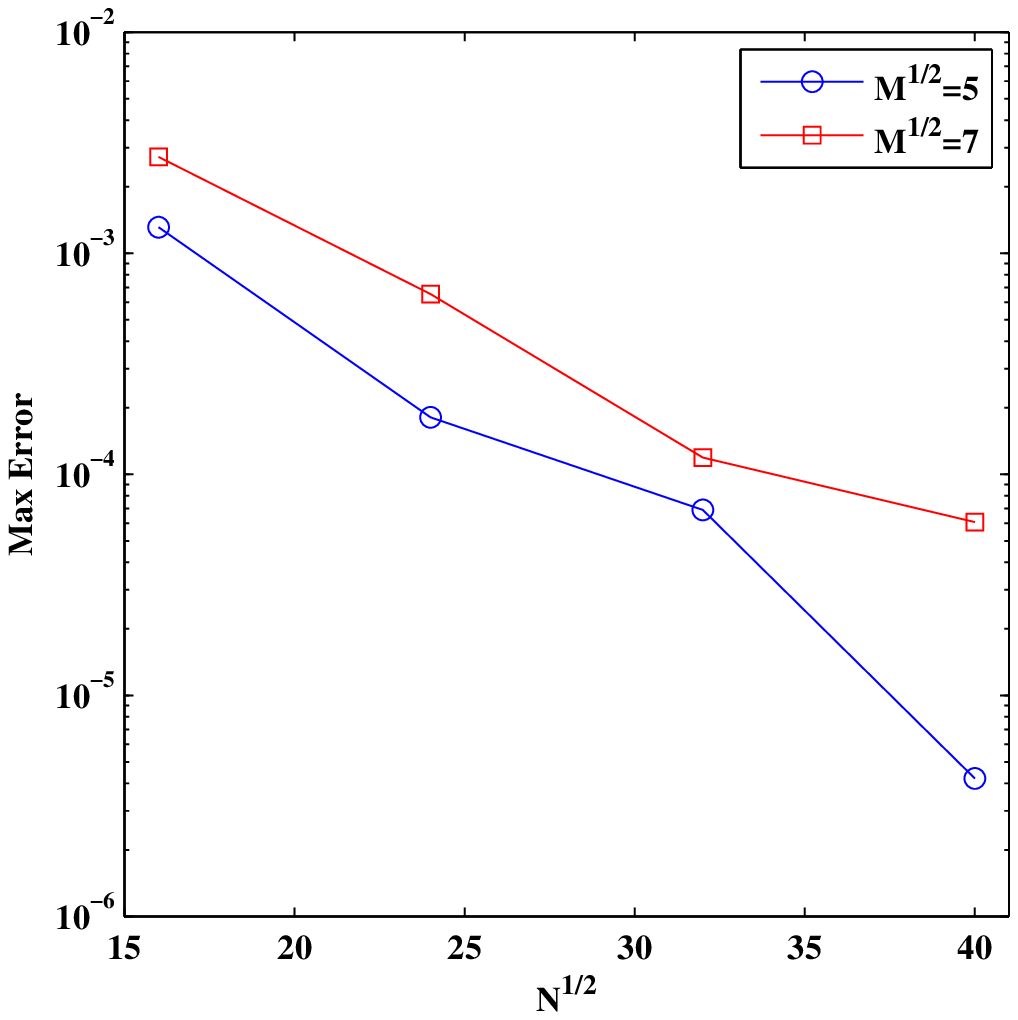}
\captionof*{figure}{Uniform points: M4, $\epsilon = 0.35$}}
\quad
\parbox{5cm}{\centering
\includegraphics[scale=0.37]{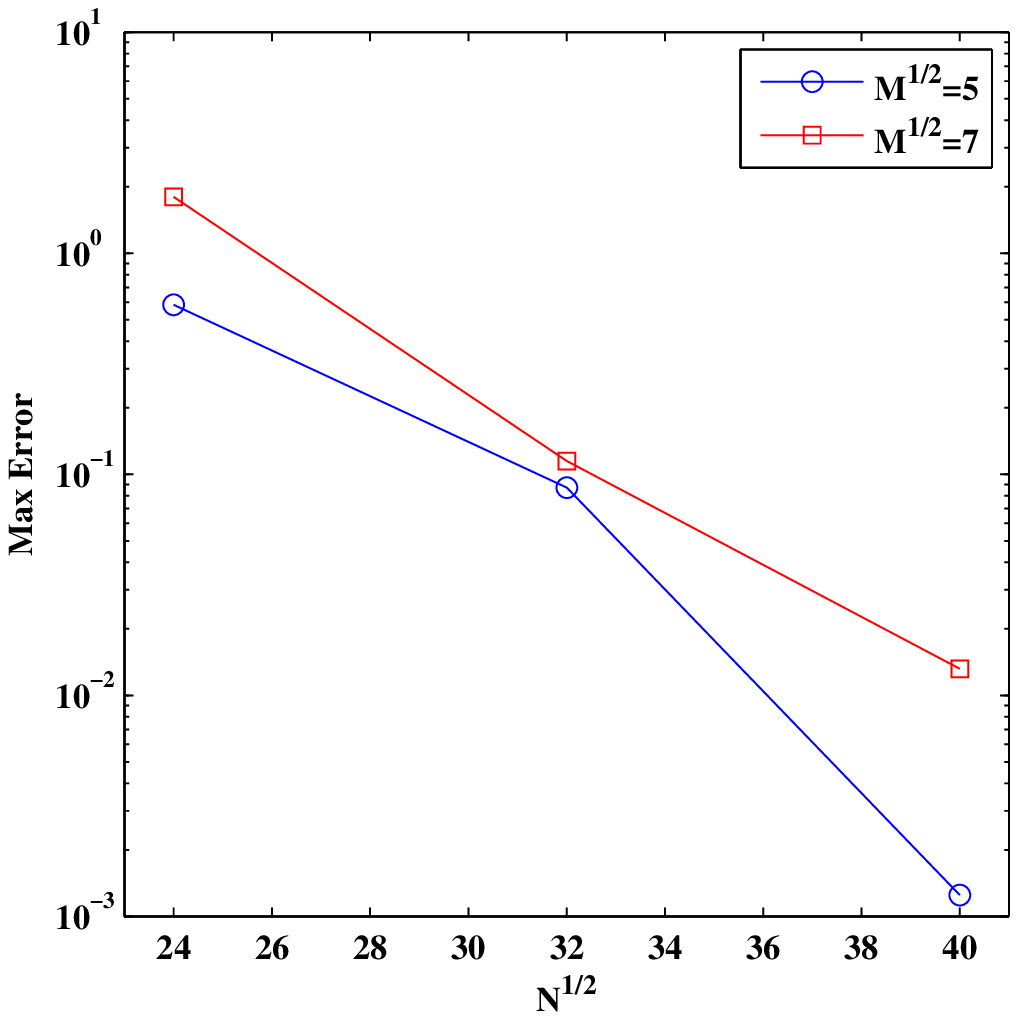}
\captionof*{figure}{Uniform points: GA, $\epsilon = 8.80$}}\\
\parbox{5cm}{\centering
\includegraphics[scale=0.37]{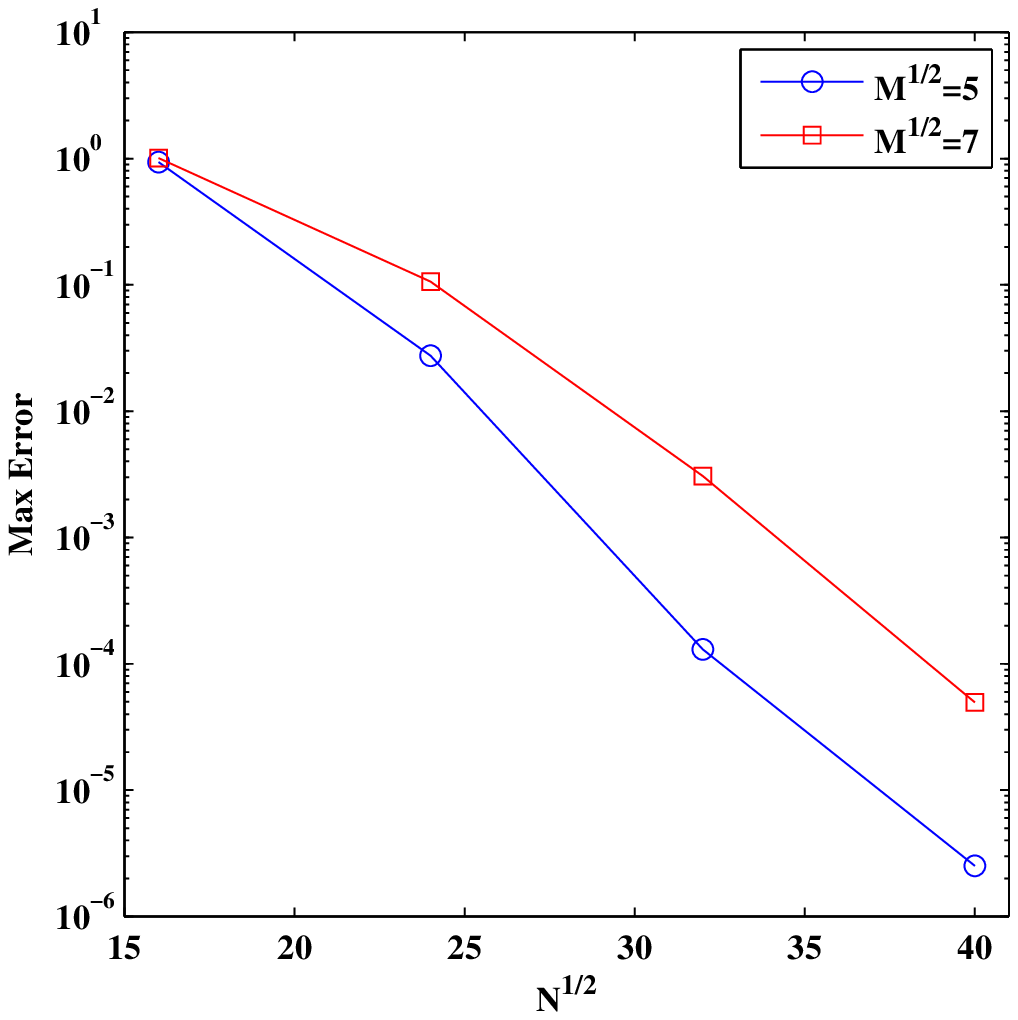}
\captionof*{figure}{Halton points: IMQ, $\epsilon = 3.05$}}
\quad
\parbox{5cm}{\centering
\includegraphics[scale=0.37]{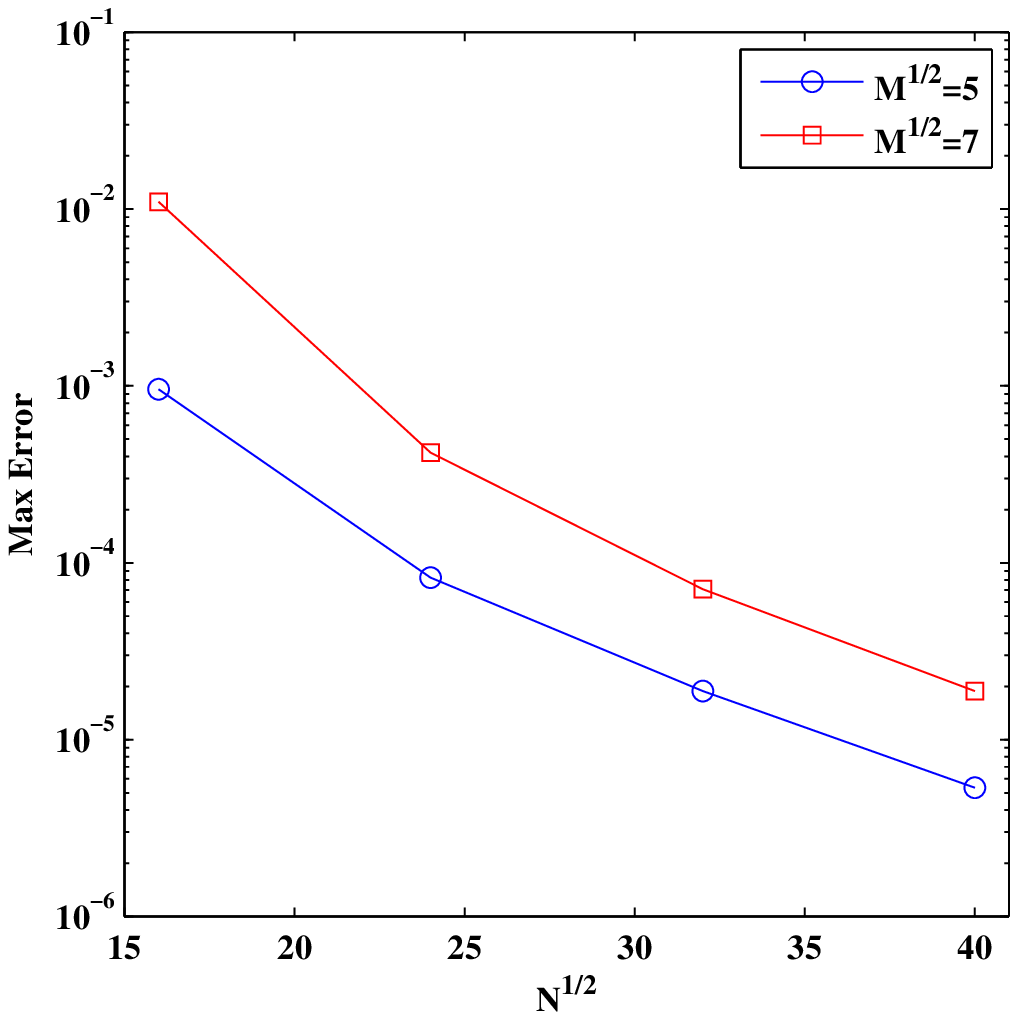}
\captionof*{figure}{Halton points: M4, $\epsilon = 0.15$}}
\quad
\parbox{5cm}{\centering
\includegraphics[scale=0.37]{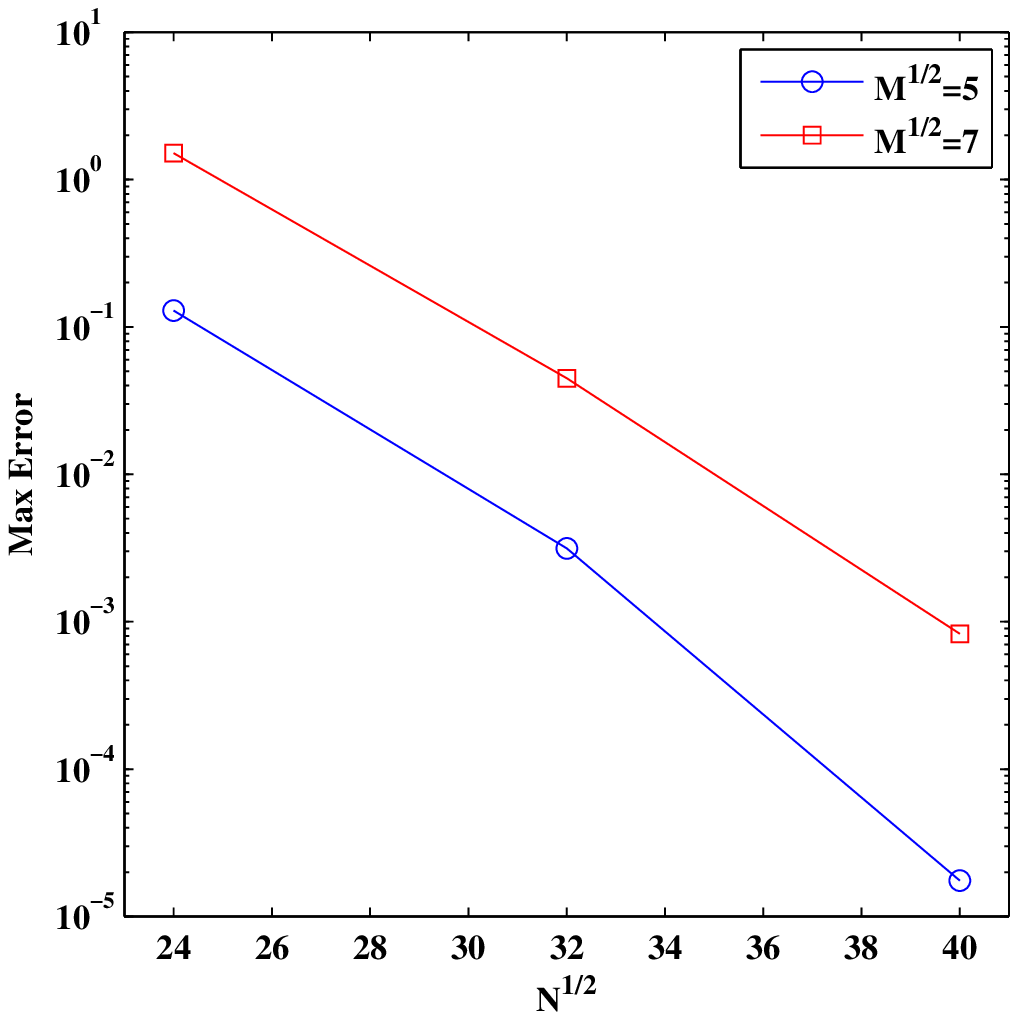}
\captionof*{figure}{Halton points: GA, $\epsilon = 7.20$}}
\caption{Error in the pseudo-papabolic equation against the problem size with respect to the number of partitions using $\delta t=0.001$ and $T=1$.}\label{f:2b}
\end{figure}

\begin{figure}
\centering
\parbox{7.5cm}{\centering
\includegraphics[scale=0.325]{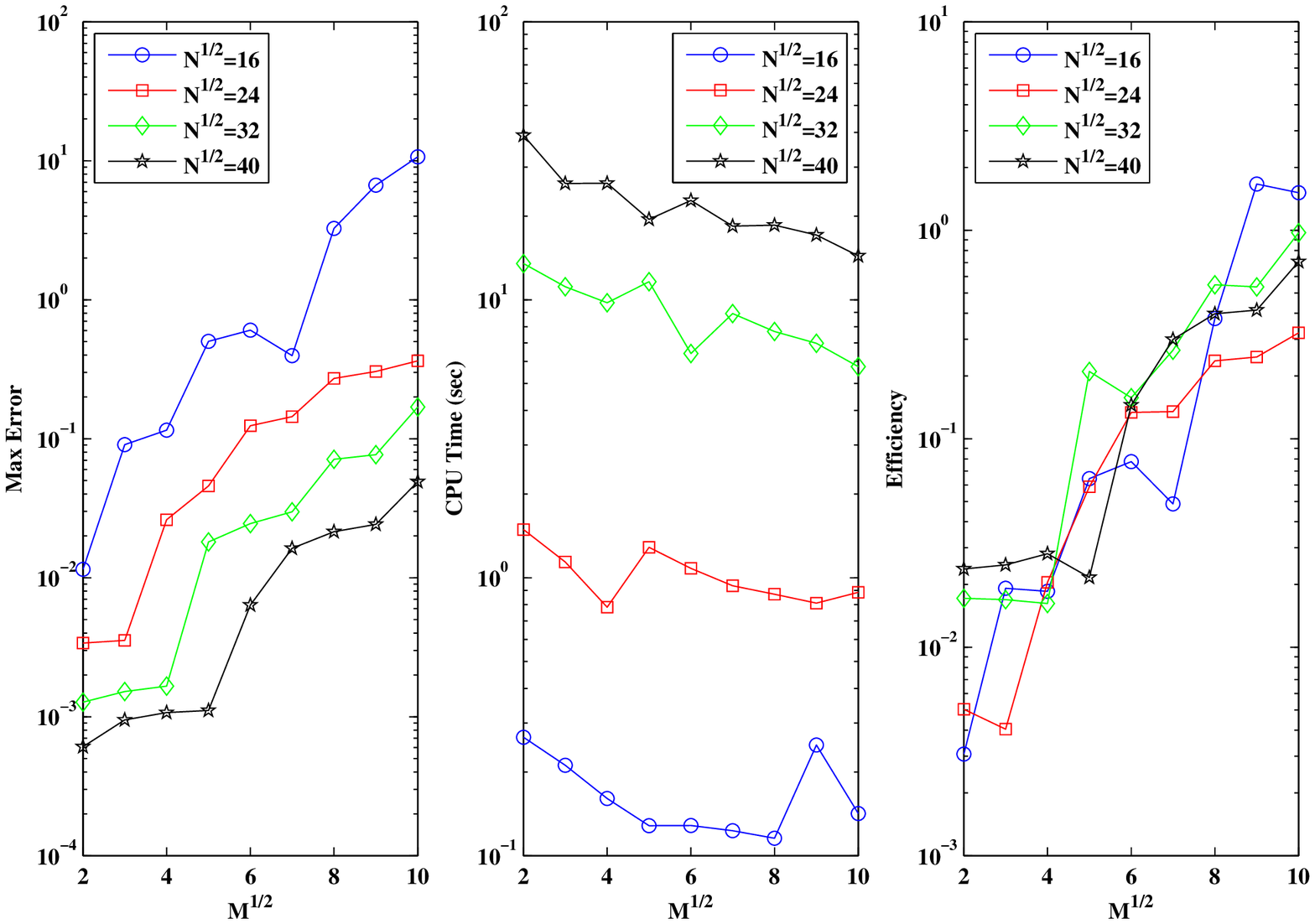}
\captionof*{figure}{Uniform points: M4, $\epsilon = 5$}}
\parbox{7.5cm}{\centering
\includegraphics[scale=0.325]{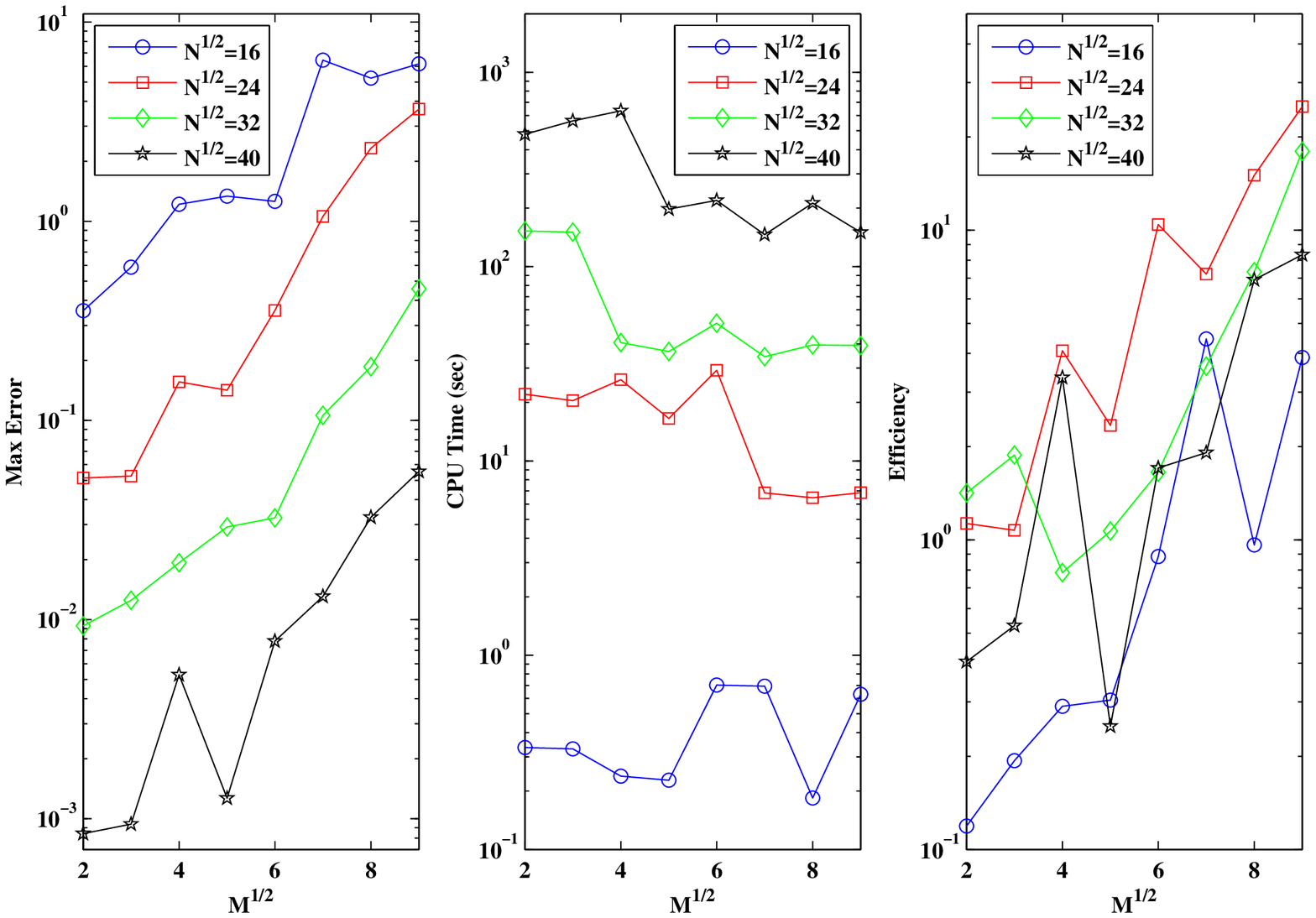}
\captionof*{figure}{Halton points: IMQ, $\epsilon = 7$}}
\caption{Left: Error in the pseudo-papabolic equation against the number of partitions in one spatial dimension. Centre: Computational time against the number of partitions in one spatial dimension. Right: Efficiency computed as product between the error and CPU time. Parameters: $\delta t=0.001$ and $T=1$.}\label{f:1b}
\end{figure}

As in Subsection \ref{res1}, in Figure \ref{f:2b} we examine how the MAE varies for Eq. \eqref{e9-1} in the RBF-PUM-FD scheme with respect to the number of points. More precisely, we analyze what happens in two cases, i.e. taking as number of subdomains (along each direction) $M^{1/2} = 5, 7$. This study involves three RBFs (IMQ, GA and M4), which are used as local approximants in the RBF-PUM-FD collocation method for uniform and Halton points. Moreover, as expected, these tests point out a similar behavior of error compared to the convection-diffusion case, thus confirming the influence of subdomain size on the accuracy of our method. It follows that covering the domain with small (large) patches results usually in less (more) accurate results even if computationally less (more) expensive. This fact appears clearly from Figure \ref{f:1b}, where by varying the number of points we show MAE (left), CPU time (centre), and efficiency (right) for a different number of partitions of $\Omega$. As earlier, for brevity we provide the results obtained in only one case for uniform and Halton points, respectively. 

\begin{figure}
\centering
\includegraphics[scale=0.6]{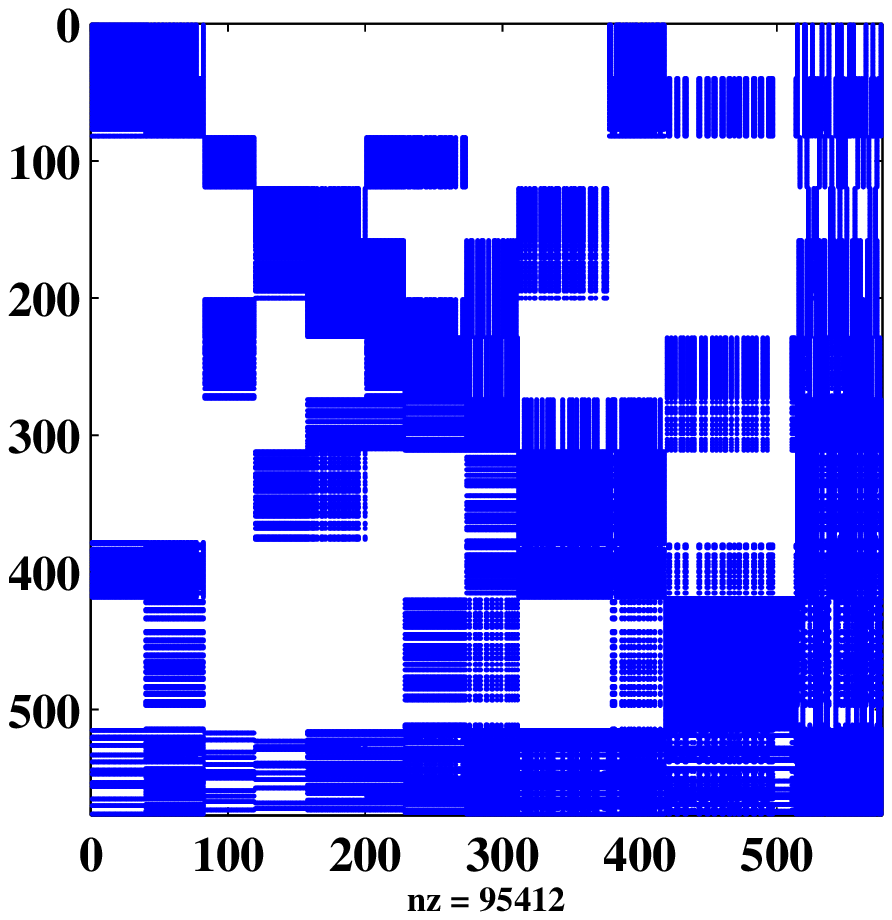}
\hskip -1.5cm
\includegraphics[scale=0.6]{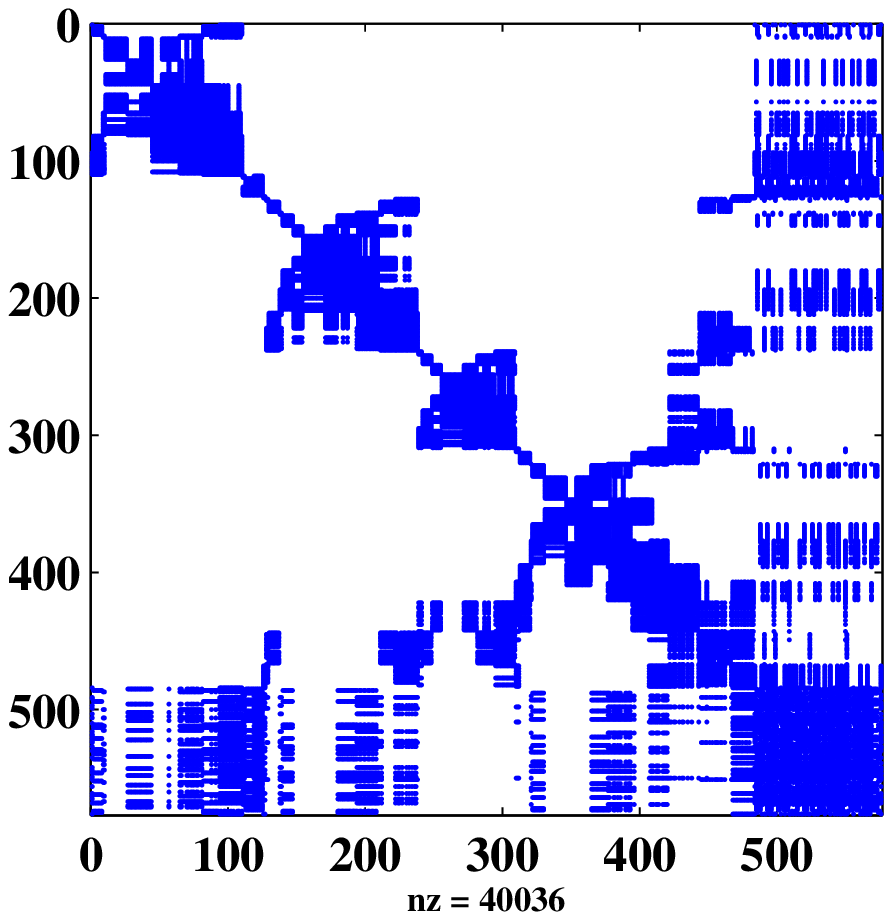}
\caption{Sparsity structure of cofficient matrix regarding the pseudo-parabolic equation with $24\times24$ Halton nodes and $3\times3$ patches (left) and $5\times5$ patches (right) for IMQ with $\epsilon=1.85$.}\label{f1a:1}
\end{figure}

Moreover, in Figure \ref{f1a:1} we show sparse structure of the matrix generated by using the RBF-PUM-FD collocation method with local IMQ approximations and two different domain partitions for the pseudo-parabolic equation. As in the convection-diffusion case, although the matrix structure is different, we can derive similar conclusions.

\begin{figure}
\centering
\includegraphics[scale=0.5]{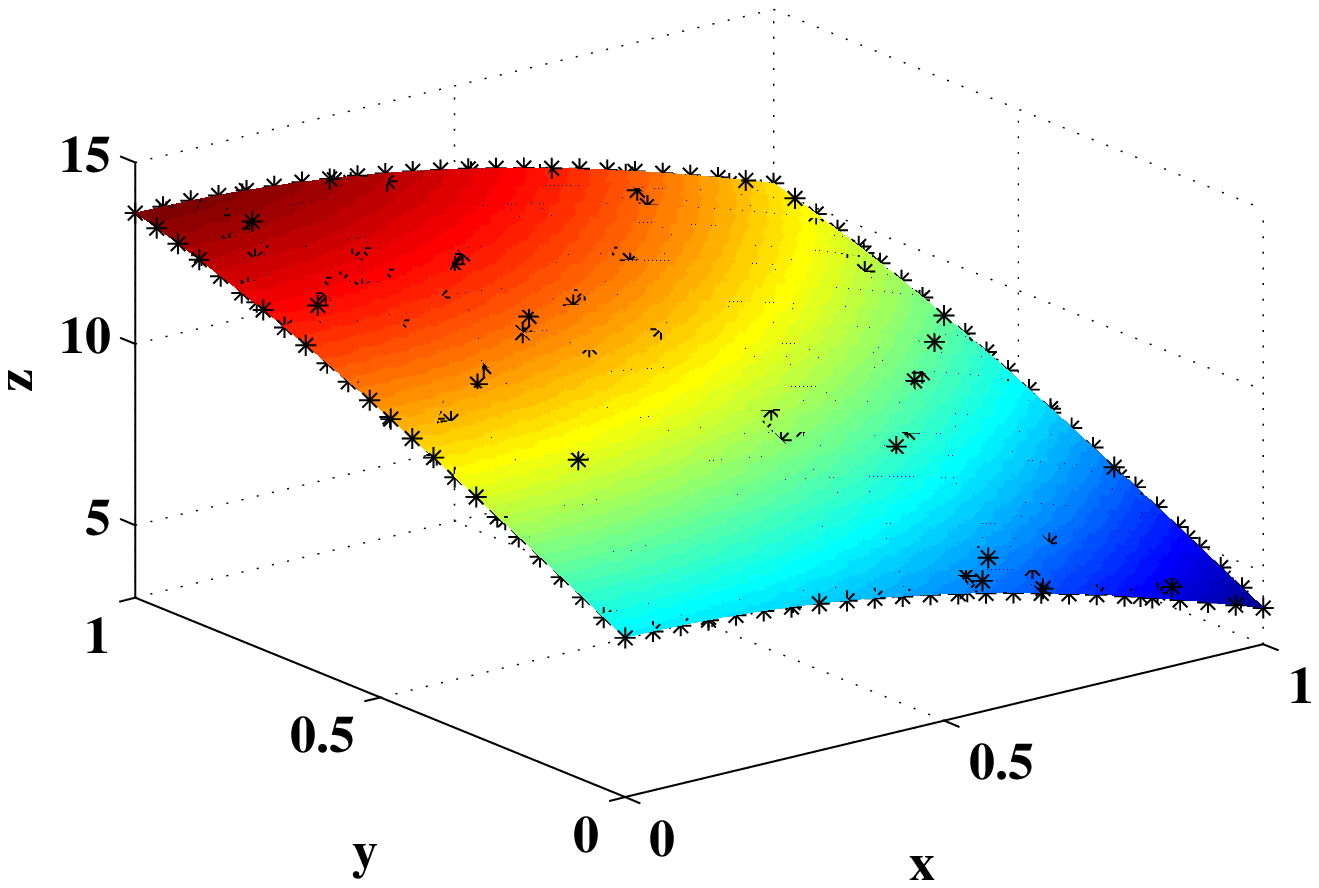}
\quad
\includegraphics[scale=0.5]{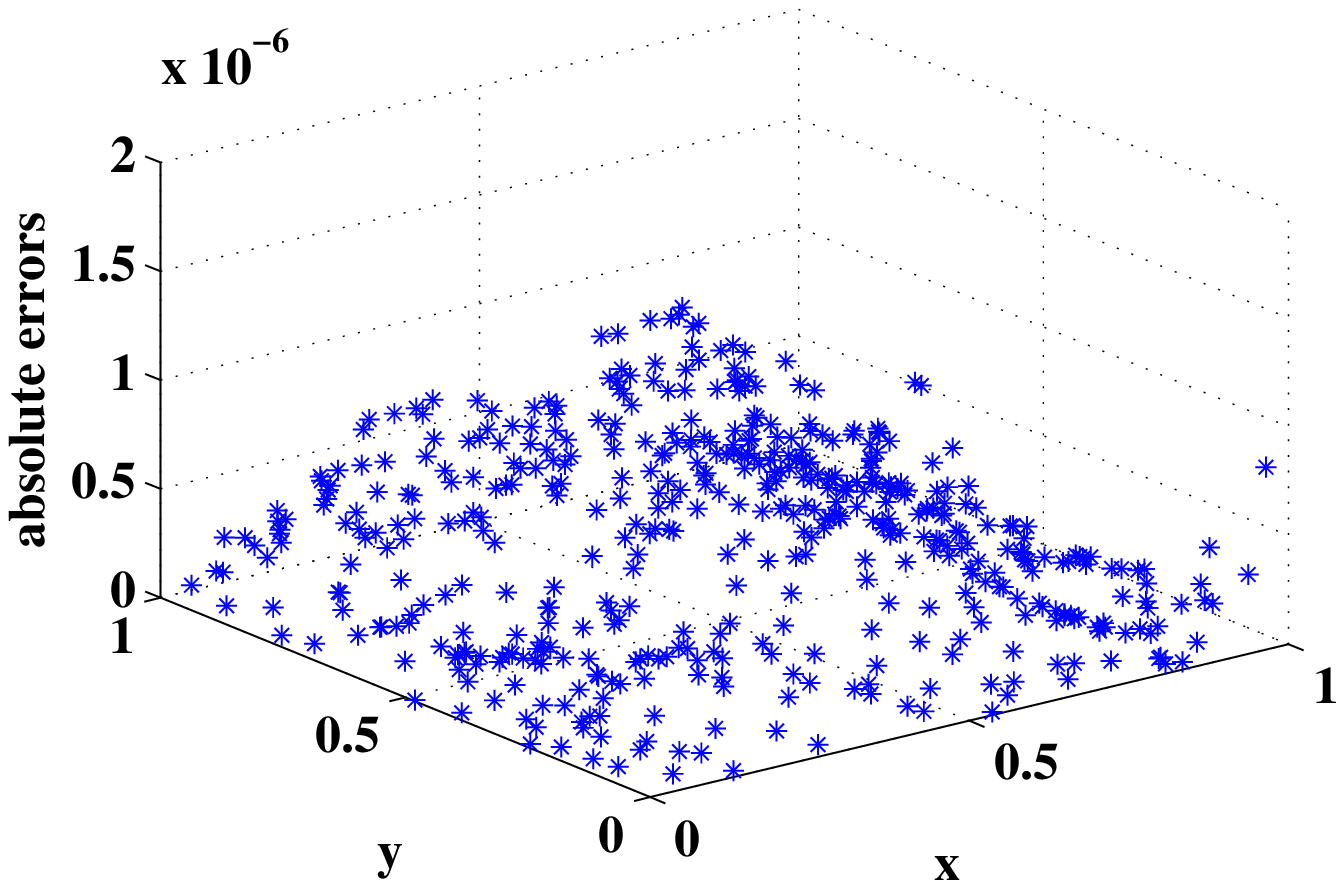}
\caption{Surface of exact solution \eqref{e9-1} and its approximation given per points (left) along with the related absolute errors (right) on $24\times24$ Halton nodes and $3\times3$ patches, using IMQ for $\epsilon = 1.80 $ with $\delta t=0.001$ and $T=1$.}
\label{f1::1-h}
\end{figure}

Then, in Tables \ref{tab:1b}-\ref{tab:1b-g-h} we compare our RBF-PUM-FD scheme with the RBF-FD method studied in \cite{esm16}. Therefore, for some sets of uniform and Halton points we show the numerical results obtained by using IMQ and GA as local RBF approximants. Specifically, we report the MAE with the associated \lq\lq optimal\rq\rq\ value of $\epsilon$, the time computed in seconds for the system solution and the CN. In addition, taking a set of Halton points and using IMQ, Figure \ref{f1::1-h} shows the surface of exact and approximate solutions (left) as well as the absolute errors (right) computed at the collocation points. As time parameters we choose again the values of $\delta t = 0.001$ and $T=1$. From such tables, we observe that accuracy of the RBF-PUM-FD method is at least comparable but often better than the one of the RBF-FD scheme. Additionally, the proposed method offers a remarkable decrease of execution time and condition number compared to the global approach.

\begin{table}
\begin{center} 
{\small
\begin{tabular}{ccccccccccc} \toprule
 & & & & & RBF-PUM-FD & $\quad$ & & & & RBF-FD \\ \cmidrule{2-6} \cmidrule{8-11}
$\sqrt{N}$ & $\sqrt{M}$& $\epsilon$ & MAE & time & CN & $\quad$ &$\epsilon$ & MAE & time & CN \\ \midrule
$16$ & $2$& $1.35$ & $5.33 \times 10^{-5}$ & $0.25$ & $6.37\times 10^{+01}$ & $\quad$ & $1.60$ & $1.49 \times 10^{-4}$ & $7.49$ & $1.20 \times 10^{+17}$\\ \midrule
$24$ & $3$& $1.70$ & $9.81 \times 10^{-6}$ & $5.48$ & $3.71\times 10^{+02}$ & $\quad$ & $2.50$ & $1.62 \times 10^{-4}$ & $22.57$ & $4.11 \times 10^{+17}$\\ \midrule
$32$ & $4$& $2.85$ & $5.29 \times 10^{-5}$ & $14.56$ & $3.22\times 10^{+02}$ & $\quad$ & $3.45$ & $1.78 \times 10^{-4}$ & $157.89$ & $1.77 \times 10^{+17}$\\ \midrule
$40$ & $5$& $3.20$ & $1.99 \times 10^{-5}$ & $31.46$ & $1.19\times 10^{+04}$ & $\quad$ & $4.15$ & $1.17 \times 10^{-4}$ & $640.12$ & $7.87 \times 10^{+17}$\\ \bottomrule
\end{tabular}
}
\caption{\small Comparison between RBF-PUM-FD and RBF-FD \cite{esm16} using IMQ for the pseudo-parabolic equation: MAE, CPU time (in seconds) and CN for uniform points with $\delta t=0.001$ and $T=1$.}
\label{tab:1b}
\end{center}
\end{table}
\begin{table}
\begin{center}
{\small
\begin{tabular}{ccccccccccc} \toprule
 & & & & & RBF-PUM-FD & $\quad$ & & & & RBF-FD \\ \cmidrule{2-6} \cmidrule{8-11}
$\sqrt{N}$ & $\sqrt{M}$& $\epsilon$ & MAE & time & CN & $\quad$ &$\epsilon$ & MAE & time & CN \\ \midrule
$16$ & $2$& $1.05$ & $9.37 \times 10^{-7}$ & $2.04$ & $8.57\times 10^{+03}$ & $\quad$ & $1.30$ & $1.52 \times 10^{-6}$ & $7.28$ & $6.93 \times 10^{+17}$\\ \midrule
$24$ & $3$& $1.80$ & $1.54 \times 10^{-6}$ & $19.17$ & $2.65\times 10^{+05}$ & $\quad$ & $2.20$ & $3.87 \times 10^{-6}$ & $42.83$ & $1.08 \times 10^{+18}$\\ \midrule
$32$ & $4$& $2.40$ & $1.28 \times 10^{-6}$ & $32.55$ & $1.91\times 10^{+06}$ & $\quad$ & $3.00$ & $6.56 \times 10^{-6}$ & $184.71$ & $1.04 \times 10^{+18}$\\ \midrule
$40$ & $5$& $3.05$ & $2.52 \times 10^{-6}$ & $81.42$ & $2.75\times 10^{+05}$ & $\quad$ & $3.75$ & $2.60 \times 10^{-6}$ & $585.13$ & $2.43 \times 10^{+18}$\\ \bottomrule
\end{tabular}
}
\caption{\small Comparison between RBF-PUM-FD and RBF-FD \cite{esm16} using IMQ for the pseudo-parabolic equation: MAE, CPU time (in seconds) and CN for Halton points with $\delta t=0.001$ and $T=1$.}
\label{tab:1b-h}
\end{center}
\end{table}
\begin{table}
\begin{center} 
{\small
\begin{tabular}{ccccccccccc} \toprule
 & & & & & RBF-PUM-FD & $\quad$ & & & & RBF-FD \\ \cmidrule{2-6} \cmidrule{8-11}
$\sqrt{N}$ & $\sqrt{M}$& $\epsilon$ & MAE & time & CN & $\quad$ &$\epsilon$ & MAE & time & CN \\ \midrule
$16$ & $2$& $3.00$ & $5.58 \times 10^{-4}$ & $0.64$ & $6.88\times 10^{+02}$ & $\quad$ & $4.00$ & $5.11 \times 10^{-3}$ & $6.13$ & $5.67 \times 10^{+17}$\\ \midrule
$24$ & $3$& $4.95$ & $8.84 \times 10^{-4}$ & $4.33$ & $1.47\times 10^{+03}$ & $\quad$ & $6.95$ & $1.40 \times 10^{-2}$ & $32.70$ & $4.13 \times 10^{+17}$\\ \midrule
$32$ & $4$& $7.60$ & $1.46 \times 10^{-3}$ & $26.68$ & $6.59\times 10^{+03}$ & $\quad$ & $9.85$ & $2.13 \times 10^{-2}$ & $182.26$ & $2.90 \times 10^{+17}$\\ \midrule
$40$ & $5$& $8.80$ & $1.25 \times 10^{-3}$ & $55.61$ & $7.60\times 10^{+04}$ & $\quad$ & $11.75$ & $1.29 \times 10^{-2}$ & $531.12$ & $5.52 \times 10^{+18}$\\ \bottomrule
\end{tabular}
}
\caption{\small Comparison between RBF-PUM-FD and RBF-FD \cite{esm16} using GA for the pseudo-parabolic equation: MAE, CPU time (in seconds) and CN for uniform points with $\delta t=0.001$ and $T=1$.}
\label{tab:1b-g}
\end{center}
\end{table}
\begin{table}
\begin{center}
{\small
\begin{tabular}{ccccccccccc} \toprule
 & & & & & RBF-PUM-FD & $\quad$ & & & & RBF-FD \\ \cmidrule{2-6} \cmidrule{8-11}
$\sqrt{N}$ & $\sqrt{M}$& $\epsilon$ & MAE & time & CN & $\quad$ &$\epsilon$ & MAE & time & CN \\ \midrule
$16$ & $2$& $2.25$ & $2.76 \times 10^{-6}$ & $1.80$ & $1.16\times 10^{+05}$ & $\quad$ & $3.35$ & $1.70 \times 10^{-5}$ & $7.21$ & $4.26 \times 10^{+17}$\\ \midrule
$24$ & $3$& $4.05$ & $1.71 \times 10^{-5}$ & $23.24$ & $2.39\times 10^{+05}$ & $\quad$ & $6.05$ & $2.20 \times 10^{-4}$ & $40.97$ & $3.52 \times 10^{+17}$\\ \midrule
$32$ & $4$& $5.50$ & $2.50 \times 10^{-5}$ & $54.84$ & $2.52\times 10^{+05}$ & $\quad$ & $8.55$ & $4.46 \times 10^{-4}$ & $158.40$ & $1.67 \times 10^{+17}$\\ \midrule
$40$ & $5$& $7.20$ & $1.75 \times 10^{-5}$ & $84.84$ & $1.54\times 10^{+06}$ & $\quad$ & $10.85$ & $3.57 \times 10^{-4}$ & $505.09$ & $1.98 \times 10^{+17}$\\ \bottomrule
\end{tabular}
}
\caption{\small Comparison between RBF-PUM-FD and RBF-FD \cite{esm16} using GA for the pseudo-parabolic equation: MAE, CPU time (in seconds) and CN for Halton points with $\delta t=0.001$ and $T=1$.}
\label{tab:1b-g-h} 
\end{center}
\end{table}
   
Finally, in Table \ref{tab:tb} we conclude this section showing the total execution times for the RBF-PUM-FD and RBF-FD methods. From these results, we can deduce considerations similar to those explained at the end of the previous subsection.

\begin{table}
\begin{center}  
{\small
\begin{tabular}{cccccccc} \toprule
 & & Uniform points & & $\quad$ & & Halton points & \\ \cmidrule{2-4} \cmidrule{6-8}
$\sqrt{N}$ & $\sqrt{M}$ & RBF-PUM-FD & RBF-FD & $\quad$ & $\sqrt{M}$ & RBF-PUM-FD & RBF-FD \\ \midrule
$16$ & $2$ & $1.38$ & $7.61$ & $\quad$ & $2$ & $4.05$ & $7.39$ \\ \midrule
$24$ & $3$ & $24.24$ & $24.88$ & $\quad$ & $3$ & $30.43$ & $42.83$ \\ \midrule
$32$ & $4$ & $30.30$ & $159.41$ & $\quad$ & $4$ & $84.15$ & $186.22$ \\ \midrule
$40$ & $5$ & $72.03$ & $643.80$ & $\quad$ & $5$ & $160.68$ & $588.68$ \\ \bottomrule
\end{tabular}
}
\caption{\small Comparison of total CPU time (in seconds) between RBF-PUM-FD and RBF-FD \cite{esm16} using IMQ for the pseudo-parabolic equation with $\delta t=0.001$ and $T=1$.} 
\label{tab:tb}
\end{center}
\end{table}


\section{Conclusions and future work}\label{sec7}

In this work, we implemented and examined a RBF-PUM-FD scheme used to solve two difficult benchmark problems such as those typically arising in various scientific and engineering applications. From this study, it emerges that the RBF-PUM-FD collocation method allows to overcome the high computational cost associated with the global RBF method, while maintaining high accuracy. Additionally, the RBF-PUM-FD technique also enables to reach a given level of accuracy with significantly less computational effort than the global RBF-FD method and with remarkably more stability than the local RBF-PUM method. The fact that RBF methods are meshfree permits an easy implementation of adaptive grids, which can be clustered around critical regions such as the strike area or the free boundary, in order to improve accuracy or reduce overall computational cost. In the case of RBF-PUM-FD, refinements can be made independently within the partitions, increasing the flexibility. However, this topic is out of the scope of the present paper. Therefore, the development of a technique for automatic adaptivity will be part of our future work, where we will also consider higher-dimensional computational problems studying and analyzing theoretically the invertibility issue of the coefficient matrix of our scheme.


\section*{Acknowledgments}
The second author acknowledges support from the Department of Mathematics \lq\lq Giuseppe Peano\rq\rq of the University of Torino via 2016-2017 project \lq\lq Multivariate approximation and efficient algorithms with applications to algebraic, differential and integral problems\rq\rq. Major parts of this research were performed while the first author visited the University of Torino supported by the Ministry of Science, Research and Technology of the Islamic Republic of Iran through a scholarship.




\end{document}